%% file: fcoho.tex
\documentclass{article}
\usepackage{varioref}
\usepackage{srcltx}
\input{xy}
\xyoption{all}

\usepackage{varioref}
\usepackage{amsfonts,
amsthm,
amsmath,latexsym,amscd,amssymb,doc}
\usepackage[latin1]{inputenc}
\usepackage{graphicx}
\usepackage[a4paper]{hyperref}

\include{header_categories} 
\include{header_sections}

  \newcommand{\gl}{\index}

  \newcommand{\Dph}[1]{\emph{#1}}
  \newcommand{\Def}[1]{\Dph{#1}\gl{#1}}

\newcommand{\mylabel}[1]{\label{#1}}

\newcommand{\status}[1]{}

\newcommand{\OF}{{\mathcal{O}_F}}
\newcommand{\OFpp}{{\mathcal{O}_{F_\pp}}}

\newcommand{\dual}{^\vee}

\newcommand{\C}{\mathcal{C}}

\newcommand{\R}{\mathrm{R}}
\newcommand{\RG}{{\R} {\Gamma}}

\newcommand{\lr}{{\longrightarrow}}
\renewcommand{\r}{\rightarrow}

\renewcommand{\t}{{\otimes}}

\newcommand{\A}{\mathbb{A}^1}
\renewcommand{\P}[1][1]{\mathbb P^{#1}}
\newcommand{\N}{\mathbb{N}}
\newcommand{\Q}{\mathbb{Q}}

\newcommand{\Ql}{{\mathbb{Q}_{{\ell}}}}
\newcommand{\Zl}{{\mathbb{Z}_{{\ell}}}}

\newcommand{\CC}{\mathbb{C}}

\newcommand{\Z}{\mathbb{Z}}

\newcommand{\pp}{\mathfrak{p}}
\newcommand{\Fpp}{{\mathbb{F}_{\pp}}}

\newcommand{\Fq}{{\mathbb{F}_q}}

\renewcommand{\H}{\mathrm{H}}

\newcommand{\p}{{{}^\mathrm{p}}}
\newcommand{\pH}{{\p\H}}

\newcommand{\x}{{\times}}
\newcommand{\ol}[1]{{\overline{#1}}} 
\newcommand{\one}{\mathbf{1}}

\newcommand{\Beweis}{{\normalfont} \textbf{Proof}}

\newcommand{\lc}{\textit{loc.~cit.}}
\newcommand{\oc}{\textit{op.~cit.}}

\newcommand{\Spec}{\mathrm{Spec}\text{ }}
\newcommand{\SpecOF}{{\Spec} {\OF}}
\newcommand{\crys}{\operatorname{crys}}
\newcommand{\SpecZ}{{\Spec} {\Z}}

\newcommand{\SpecF}{{\Spec} {F}}

\newcommand{\SpecFpp}{{\Spec} {\Fpp}}

\newcommand{\eff}{\mathrm{eff}}

\newcommand{\red}{\mathrm{red}}

\newcommand{\gr}{\operatorname{gr}} 

\newcommand{\h}{\mathrm{h}}

\newcommand{\wt}{\operatorname{wt}} 
\newcommand{\num}{\mathrm{num}} 
\renewcommand{\hom}{\mathrm{hom}} 
\newcommand{\rat}{\mathrm{rat}} 
\newcommand{\Hom}{\mathrm{Hom}}
\newcommand{\IHom}{\underline{\Hom}}
\newcommand{\id}{\mathrm{id}}
\renewcommand{\gr}{\operatorname{gr}} 
\newcommand{\RHom}{{\R}{\Hom}}

\newcommand{\im}{\operatorname{im}} 
\renewcommand{\Im}{\mathrm{im}}

\newcommand{\M}{\operatorname{M}} 
\renewcommand{\h}{\mathrm{h}} 
\newcommand{\Mc}{\M_{\mathrm c}} 

\newcommand{\Gal}{\mathrm{Gal}}
\newcommand{\CH}{\mathrm{CH}}

\newcommand{\bound}{\mathrm{b}} 

\newcommand{\Ext}{\operatorname{Ext}}

\newcommand{\loc}[2]{[}

\newcommand{\pr}{\begin{proof}[\Beweis: ]}
\newcommand{\pf}{\pr}
\newcommand{\xpf}{\end{proof}}
\newcommand{\xproof}{\end{proof}}

\makeindex

\makeatletter
\@addtoreset{Defi}{chapter}
\@addtoreset{Defi}{section}
\makeatother

\begin{document}

\title{$f$-cohomology and motives over number rings}
\author{Jakob Scholbach \footnote{Universit{\"a}t M{\"u}nster, Mathematisches Institut, Einsteinstr. 62, D-48149 M{\"u}nster, Germany, 
\href{mailto:jakob.scholbach@uni-muenster.de}{jakob.scholbach@uni-muenster.de}
}}
\maketitle

\begin{abstract}
This paper is concerned with an interpretation of $f$-cohomology, a modification of motivic cohomology of motives over number fields, in terms of motives over number rings. Under standard assumptions on mixed motives over finite fields, number fields and number rings, we show that the two extant definitions of $f$-cohomology of mixed motives $M_\eta$ over a number field $F$---one via ramification conditions on $\ell$-adic realizations, another one via the $K$-theory of proper regular models---both agree with motivic cohomology of $\eta_{!*} M_\eta[1]$. Here $\eta_{!*}$ is constructed by a limiting process in terms of intermediate extension functors $j_{!*}$ defined in analogy to perverse sheaves. 
\end{abstract}

The aim of this paper is to give an interpretation of $f$-cohomology in terms of motives over number rings. The notion of $f$-cohomology goes back to Beilinson who used it to formulate a conjecture about special $L$-values \cite{Beilinson:Higher,Beilinson:Notes}. The most classical example is what is now called $\H^1_f(F, \one(1))$, $f$-cohomology of $\one(1)$, the motive of a number field $F$, twisted by one. This group is $\mathcal O_F^\x \t_\Z \Q$, as opposed to the full motivic cohomology  $\H^1(F, \one(1)) = F^\x \t \Q$. Together with the Dirichlet regulator, it explains the residue of the Dedekind zeta function $\zeta_F(s)$ at $s=1$. This idea has been generalized in many steps and many ways, for example to the notion of Selmer complexes \cite{Nekovar:Selmer}. This work is concerned with the $f$-cohomology of a mixed motive $M_\eta$ over $F$. There are two independent yet conjecturally equivalent ways to define $\H^1_f(F, M_\eta) \subset \H^1(F, M_\eta)$. We interpret the two definitions of $f$-cohomology as motivic cohomology of suitable motives over $\OF$. This idea is due to Huber. 

There are two approaches to $\H^1_f(M_\eta)$. The first is due to Beilinson \cite[Remark 4.0.1.b]{Beilinson:Height}, Bloch and Kato \cite[Conj. 5.3.]{BlochKato} and Fontaine \cite{Fontaine:Valeurs,FPR}. It is given by picking elements in motivic cohomology acted on by the local Galois groups in a prescribed way (\refde{H1flocal}, \refde{H1fglobal}, \refde{H1fEllAdic}). 
The second definition of $\H^1_f(M_\eta)$, due to Beilinson \cite[Section 8]{Beilinson:Notes}, applies to $M_\eta = \h^{i-1}(X_\eta)(n)$, with $X_\eta$ smooth and projective over $F$, $i - 2n < 0$. It is given by the image of $K$-theory of a regular proper model $X$ of $X_\eta$ (\refde{H1fImage}). 
Such a model may not exist, but there is a unique meaningful extension of this definition to all Chow motives over $F$ due to Scholl \cite{Scholl:Integral}. 

Our main results (Theorems \ref{theo_charH1fmotivic}, \ref{theo_summaryH1f}, \ref{theo_Scholl}) show that both definitions of $\H^1_f(M_\eta)$ agree with $\H^0(\eta_{!*} \h^{i-1}(X_\eta, n)[1])$. Here $\eta_{!*}$ is a functor that attaches to any suitable mixed motive over $F$ one over $\OF$. It is defined by a limiting process using the intermediate extension $j_{!*}$ familiar from perverse sheaves \cite{BBD} along open immersions $j: U \r \SpecOF$. Even to formulate such a definition, one has to rely on profound conjectures, namely the existence of mixed motives over (open subschemes of) $\SpecOF$. The proof of the main theorems also requires us to assume a number of properties related to weights of motives.

We point out that previously Jannsen and Scholl have shown the agreement of these two notions (in the case $M_\eta = h^{i}(X_\eta)(n)$, $X_\eta / F$ smooth and proper) under weaker hypotheses than the ones considered here \cite{Scholl:IntegralII}. Also Scholl unconditionally proved the agreement for products of smooth projective curves over $F$ (\oc). Our motivation for this studying and employing this stronger set of assumptions about motives lies in an application to special $L$-values conjectures \cite{Scholbach:specialL}. Very briefly, Beilinson's conjecture concerning special $L$-values for mixed motives $M_\eta$ over $\Q$ has $f$-cohomology as motivic input. $L$-functions of such motives can be generalized to motives over $\Z$ such that the classical $L$-function of $M_\eta$ agrees with the $L$-function (over $\Z$) of $\eta_{!*} M_\eta[1]$. Thereby the $L$-function and the motivic data in Beilinson's conjecture belong to the same motive over $\Z$, thus giving content to a more general conjecture about special $L$-values for motives over $\Z$. In this light it is noteworthy that $\H^0(\eta_{!*} \h^{2n-1}(X_\eta, n)[1])$ identifies with the group that occurs in the part of Beilinson's conjecture that describes special values at the central point. 

The contents of the paper are as follows: \refsect{axiomaticdescription} is the basis of the remainder; it lists a number of axioms on triangulated categories of motives. Such categories $\DMgm(S)$ have been constructed by Voevodsky \cite{Voevodsky:TCM} and Hanamura \cite{Hanamura:Motives1} (over fields) and Levine \cite{Levine:MixedMotives} (over bases $S$ over a field). The various approaches are known to be (anti-)equivalent, at least for rational coefficients \cite[Section VI.2.5]{Levine:MixedMotives}, \cite[Section 4]{Bondarko:Differential}. Over more general bases $S$, the category $\DM(S)$ has been constructed by Ivorra \cite{Ivorra:Realisation} and Cisinski and D\'eglise \cite{CisinskiDeglise:Triangulated}. We sum up the properties of this construction by specifying a number of axioms concerning triangulated categories of motives that will be used in the sequel. They are concerned with the ``core'' behavior of $\DM(S)$, that is: functoriality, compacity, the monoidal structure and the relation to algebraic $K$-theory, as well as localization, purity, base-change and resolution of singularities. We work with motives with rational coefficients only, since this is sufficient for all our purposes. We use a contravariant notation for motives, that is to say the functor that maps any scheme $X$ to its motive $\M(X)$ shall be contravariant. This is in line with most pre-Voevodsky papers.

\refsect{realizations} is a very brief reminder on realizations. The existence of various realizations, due to Huber and Ivorra \cite{Huber:Mixed, Huber:Realization, Ivorra:Realisation}, is pinning down the intuition that motives should be universal among (reasonable) cohomology theories. 

After \refsect{perverse}, a brief intermezzo on perverse $\ell$-adic sheaves over $\OF$, \refsect{conjprop} spells out a number of conjectural properties (also called axioms in the sequel) of $\DMgm(S)$, where $S$ is either a finite field $\Fpp$, a number field $F$ or a number ring $\OF$. The first group of these properties centers around the existence of a category of mixed motives $\MM(S)$, which is to be the heart of the so-called motivic $t$-structure. The link between mixed motives over $\OF$ and $\Fpp$ or $F$ is axiomatized by mimicking the exactness properties familiar from perverse sheaves (\refax{exactness}). A key requirement on mixed motives is that the realization functors on motives should be exact (\refax{tstructureRealizations}). For the $\ell$-adic realization over $\SpecOF[1/\ell]$, this requires a notion of perverse sheaves over that base (\refsect{perverse}). Another important conjectural facet of mixed motives are weights. Weights are an additional structure encountered in both Hodge structures and $\ell$-adic cohomology of algebraic varieties over finite fields, both due to Deligne \cite{Deligne:Hodge3, Deligne:Weil2}. They are important in that morphisms between Hodge structures or $\ell$-adic cohomology groups are known to be strictly compatible with weights, moreover, they are respected to a certain extent by smooth maps and proper maps. It is commonly assumed that this should be the case for mixed motives, too. We show in a separate work that the $t$-structure axioms and the needed weight properties hold in the triangulated subcategory $\DATM(\OF) \subset \DMgm(\OF)$ of Artin-Tate motives (as far as they are applicable) \cite{Scholbach:Mixed}.

The remaining two sections assume the validity of the axiomatic framework set up so far. The first key notion in \refsect{motivesOF} is the intermediate extension $j_{!*} M$ of a mixed motive $M$ along some open embedding $j$ inside $\SpecOF$. This is done as in the case of perverse sheaves, due to Beilinson, Bernstein and Deligne \cite{BBD}. Quite generally, much of this paper is built on the idea that the abstract properties of mixed perverse sheaves (should) give a good model for mixed motives over number rings. Next we develop a notion of smooth motives, which is an analog of lisse \'etale sheaves. This is needed to use a limiting technique to get the extension functor $\eta_{!*}$ that extends motives over $F$ to ones over $\OF$. Finally, we apply the axiom on the exactness of $\ell$-adic realization to show that intermediate extensions commute with the realization functors. This will be a stepstone in a separate work on $L$-functions of motives \cite{Scholbach:specialL}.  

\refsect{fcohomology} gives the comparison theorems on $f$-cohomology mentioned above. The two definitions of $f$-cohomology being quite different, the proofs of the comparison statements are different, too: the first is essentially based on the Hochschild-Serre spectral sequence. The crystalline case of that definition of $f$-cohomology is disregarded throughout. The second proof is a purely formal, if occasionally intricate bookkeeping of cohomological degrees and weights. 


The problem of finding a motivic interpretation of terms such as $\H^1_f(M_\eta)$ underlying the formulation of Beilinson's conjecture has been studied by Scholl \cite{Scholl:Remarks, Scholl:Integral, Scholl:IntegralII}, who develops an abelian category $\MM(F / \OF)$ of mixed motives over $\OF$ by taking mixed motives over $F$, and imposing additional non-ramification conditions. Conjecturally, the group $\Ext^i_{\MM(F / \OF)} (\one, \h^i(X_\eta, n))$ for $X_\eta / F$ smooth and projective, $i = 0$, $1$, agrees with what amounts to $\H^{i-1}(\eta_{!*} \h^{2n-1}(X_\eta, n)[1])$. 

No originality is claimed for Sections \ref{sect_axiomaticdescription}, \ref{sect_realizations}, and \ref{sect_conjprop}, except perhaps for the formulation of the relation of mixed motives over $\OF$ and $F$ and the residue fields $\Fpp$, which however is a natural and immediate translation of the theory of perverse sheaves. I would like to thank Denis-Charles Cisinski and Fr\'ed\'eric D\'eglise for communicating to me their work on $\DMgm(S)$ over general bases \cite{CisinskiDeglise:Triangulated} and Baptiste Morin for explaining me a point in \'etale cohomology. Most of all, I gratefully acknowledge Annette Huber's advice in writing my thesis, of which this paper is a part.

\section{Geometric motives} \mylabel{sect_axiomaticdescription}

Throughout this paper, $F$ is a number field, $\OF$ its ring of integers, $\pp$ stands for a place of $F$. For finite places, the residue field is denoted $\Fpp$. By scheme we mean a Noetherian separated scheme. Actually, it suffices to think of schemes of finite type over one of the rings just mentioned. In this section $S$ denotes a fixed base scheme.

This section is setting up a number of axioms describing a triangulated category $\DMgm(S)$ of geometric motives over $S$. They will be used throughout this work. As pointed out in the introduction, the material of this section is due to Cisinski and D\'eglise \cite{CisinskiDeglise:Triangulated}, who build such a category of motives using Ayoub's base change formalism \cite{Ayoub:Six1}. 

\axio \mylabel{axio_motiviccomplexes} (Motivic complexes and functoriality)
\begin{itemize}
\item 
There is a triangulated $\Q$-linear category $\DM(S)$\index{DM@$\DM(-)$}. It is called category of \Def{motivic complexes} over $S$ (with rational coefficients). It has all limits and colimits. 
\item (Tensor structure)
The category $\DM(S)$ is a triangulated symmetric monoidal category (see e.g. \cite[Part 2, II.2.1.3]{Levine:MixedMotives}). Tensor products commute with direct sums. The unit of the tensor structure is denoted $\one_S$ or $\one$\index{1@$\one$}. Also, there are internal $\Hom$-objects in $\DM$, denoted $\IHom$. The \Def{dual} $M\dual$ of an object $M \in \DM(S)$ is defined by $M\dual := \IHom(M, \one)$.
\item
For any map $f: X \r Y$ of schemes, there are pairs of adjoint functors
\eqn \mylabel{eqn_adjstar}
f^*: \DM(Y) \leftrightarrows \DM(X) : f_*
\xeqn
such that $f^* \one_Y = \one_X$ and, if $f$ is quasi-projective,
$$f_!: \DM(X) \leftrightarrows \DM(Y): f^!.$$
\end{itemize}
\xaxio

The existence of $f_!$ and $f^!$ is restricted to quasi-projective maps since the abstract construction of these functors in Ayoub's work \cite[Section 1.6.5]{Ayoub:Six1}, on which Cisinski's and D\'eglise's construction of motives over general bases \cite{CisinskiDeglise:Triangulated} relies, has a similar restriction.

Recall that an object $X$ in a triangulated category $\mathcal T$ closed under arbitrary direct sums is \Def{compact} if $\Hom(X, -)$ commutes with direct sums. The subcategory of $\mathcal T$ of compact objects is triangulated and closed under direct summands (a.k.a. a \Def{thick subcategory}) \cite[Lemma 4.2.4]{Neeman:TC}. The category $\mathcal T$ is called \Def{compactly generated} if the smallest triangulated subcategory closed under arbitrary sums containing the compact objects is the whole category $\mathcal T$.

\axio \mylabel{axio_compact} (Compact objects)
The motive $\one \in \DM(S)$ is compact. 
The functors $f^*$ and $f^!$, whenever defined, and $\t$ and $\IHom$ preserve compact objects. The same is true for $f_*$ and $f_!$ if $f$ is of finite type. 
The canonical map $M \r (M\dual)\dual$ is an isomorphism for any compact object $M$.
\xaxio

\defi \mylabel{defi_motive}
The subcategory of compact objects of $\DM(S)$ is denoted $\DMgm(S)$\index{DMgm@$\DMgm(-)$} and called the category of \Def{geometric motives} over $S$. 

For any map $f: X \r S$ of finite type, the object $\M_S(X) := \M(X) := f_* f^* \one \in \DMgm(S)$ is called the \Def{motive} of $X$ over $S$. By adjunction, $\M$ is a contravariant functor from schemes of finite type over $S$ to $\DMgm(S)$. 
For any quasi-projective $f: X \r S$, the \emph{motive with compact support}\index{motive!with compact support} of $X$, $\Mc(X)$, is defined as $f_! f^* \one \in \DMgm(S)$. 

The smallest thick subcategory of $\DM(S)$ containing the image of $\M$ is denoted $\DMgmeff(S)$ and called the category of \Def{effective geometric motives}. The closure of that subcategory under all direct sums is called the category of \Def{effective motives}, $\DMeff(S)$.

\xdefi

\axio \mylabel{axio_monoidal} (Tensor product vs. fiber product) The functor $\M$ is an additive tensor functor, i.e., maps disjoint unions of schemes over $S$ to direct sums and fiber products of schemes over $S$ to tensor products in $\DMgm(S)$. 
\xaxio

\axio \mylabel{axio_compactlygenerated} (Compact generation)
The categories $\DM(S)$ and $\DMeff(S)$ are compactly generated. 
\xaxio

The category $\DM(S)$, being closed under countable direct sums is pseudo-abelian \cite[Lemma II.2.2.4.8.1]{Levine:MixedMotives}, i.e., it contains kernels of projectors. In particular, the projector $\M(\P_S) \r \M(S) \r \M(\P_S)$ has a kernel $K$ (the first map is induced by the projection onto the base, the second map stems from the rational point $0 \in \P_S$). The object 
$$\one(-1) := K[2],$$ 
is called \Def{Tate object} or Tate motive. The resulting decomposition $\M(\P_S) = \one \oplus \one(-1)[-2]$ implies $\one(-1) \in \DMgmeff(S)$.

\axio \mylabel{axio_TateCancel} (Cancellation and Effectivity)
%
In $\DMgm(S)$ (and thus in $\DM(S)$), the Tate object $\one(-1)$ has a tensor-inverse denoted $\one(1)$\index{1(1)@$\one(1)$}. For any $M \in \DM(S)$, $n \in \Z$, set $M(n) := M \t \one(1)^{\t n}$. Then there is a canonical isomorphism called \Def{cancellation} isomorphism ($n \in \Z$, $M, N \in \DM(S)$):
$$\Hom_{\DM(S)}(M, N) \cong \Hom_{\DM(S)}(M(n), N(n)).$$

The smallest tensor subcategory of $\DMgm(S)$ that contains $\DMgmeff(S)$ and $\one(1)$ is $\DMgm(S)$. In other words, $\DMgm(S)$ is obtained from $\DMgmeff(S)$ by tensor-inverting $\one(-1)$.
\xaxio

\defi \mylabel{defi_motiviccoho}
Let $M$ be any geometric motive over $S$. We write $\H^i(M) := \H^i(S, M) := \Hom_{\DM(S)}(\one, M[i])$. For $M = \M(X)(n)$ for any $X$ over $S$ we also write $\H^i(X, n) := \H^i(\M(X)(n)) = \Hom_{\DMgm(S)} (\one, \M(X)(n)[i]) \stackrel{\text{\refeq{adjstar}}}= \Hom_{\DMgm(X)} (\one, \one(n)[i])$. This is called \Def{motivic cohomology} of $M$ and $X$, respectively. 
\xdefi

%

\axio \mylabel{axio_RelationKTheory} (Motivic cohomology vs. K-theory)\index{motivic cohomology!Relation to $K$-theory}
For any regular scheme $X$, there is an isomorphism $\H^i(X, n) \cong K_{2n-i}(X)^{(n)}_\Q$, where the right hand term  denotes the Adams eigenspace of algebraic $K$-theory tensored with $\Q$ \cite{Quillen:HAK}. 
\xaxio

This is a key property of motives, since algebraic $K$-theory is a universal cohomology theory in the sense that Chern characters map from algebraic $K$-theory to any other (reasonable) cohomology theory of algebraic varieties \cite{Gillet:RR}. For $S$ a perfect field, this axiom is given by \cite[Prop. 4.2.9]{Voevodsky:TCM} and its non-effective analogue. See also \cite[Theorem I.III.3.6.12.]{Levine:MixedMotives}. 

Recall Grothendieck's category of pure motives $\PureMot_\sim(K)$ with respect to an adequate equivalence relation $\sim$, see e.g.\ \cite[Section 4]{Andre:Motifs}. For rational equivalence they are also called Chow motives, since, for any smooth projective variety $X$ over a field $K$,
\eqn \mylabel{eqn_ChowM}
\Hom_{\PureMot_\rat(K)}(\one(-n), h(X)) = \CH^n(X),
\xeqn
where $h(X)$ denotes the Chow motive of $X$ and the right hand term is the \Def{Chow group} of cycles of codimension $n$ in $X$. This way, the above axiom models the fact \cite[2.1.4]{Voevodsky:TCM} that Chow motives are a full subcategory of $\DMgm(K)$. Under the embedding $\PureMot_\rat(K) \subset \DMgm(K)$, $h(X, n)$ maps to $\M(X)(n)[2n]$.

\bem
We do not need to assume \emph{expressis verbis} homotopy invariance (i.e., $\one \stackrel \cong \r pr_* pr^* \one \in \DMgm(S)$ for $pr: S \x \A \r S$) nor the projective bundle formula \cite[Prop. 3.5.1]{Voevodsky:TCM}. (Note, however, that $K'$-theory does have such properties.) 
\xbem

\axio \mylabel{axio_localization} (Localization)
Let $i: Z \r S$ be any closed immersion and $j: V \r S$ the open complement. The adjointness maps give rise to the following distinguished triangles in $\DM(S)$: 
$$j_! j^* \r \id \r i_* i^*,$$
$$i_* i^! \r \id \r j_* j^*.$$
(In particular, $f_* f^* \cong \id$, where $f: X_\red \r X$ denotes the canonical map of the reduced subscheme structure.) In addition, one has $j^* j_* = \id$ and $i^* i_* = \id$, equivalently $j^* i_* = i^* j_! = 0$.
\xaxio


\axio \mylabel{axio_purity} (Purity and base change)
\begin{itemize}
\item
For any quasi-projective map $f$, there is a functorial transformation of functors $f_! \r f_*$. It is an isomorphism if $f$ is projective.
\item
(\Def{Relative purity}): If $f$ is quasi-projective and smooth of constant relative dimension $d$, there is a functorial (in $f$) isomorphism $f^! \cong f^*(d)[2d]$. 
\item
(\Def{Absolute purity}): If $i : Z \r U$ is a closed immersion of codimension $c$ of two regular schemes $Z$ and $U$, there is a natural isomorphism $i^! \one \cong \one(-c)[-2c]$. 
\item
(\Def{Base change}): For any two quasi-projective maps $f$ and $g$ let $f'$ and $g'$ denote the pullback maps:
\eqn \label{eqn_basechange}
\xymatrix
{
X' \x_X Y \ar[r]^{g'} \ar[d]^{f'} &
Y \ar[d]^f \\
X' \ar[r]^g &
X
}
\xeqn
Then there is canonical isomorphism of functors 
$$f^* g_! \stackrel{\cong}{\lr} g'_! f'^*.$$ 
\end{itemize}
\xaxio

This axiom is proven by Cisinski \& D\'eglise using Ayoub's general base change formalism. See in particular \cite[1.4.11, 12]{Ayoub:Six1} for the construction of the base change map. See also \cite[Theorem I.I.2.4.9]{Levine:MixedMotives} for a similar statement in Levine's category of motives.


  

\defi \mylabel{defi_Verdierdual}
Let $f: S \r \SpecZ$ be the structural map. Assume $f$ is quasi-projective. Then $D(M) := \IHom(M, f^! \one(1)[2])$ is called \Def{Verdier dual} of $M$. 
\xdefi

By the preceding axioms, $D$ induces a contravariant endofunctor of $\DMgm(S)$. The shift and twist in the definition is motivated as follows: given some complex analytic space $X$, the Verdier dual of a sheaf $\mathcal F$ on $X$ is defined by
$$D(\mathcal F) := \underline{\RHom}_{\D(\Shv{}{X})}(\mathcal F,  f^! \Z),$$
where $f$ denotes the projection to a point, see e.g.\ \cite[Ch. VI]{Iversen:Cohomology}. When $X$ is smooth of dimension $d$, one has $f^! \Z = f^* \Z (d)[2d] = \Z(d)[2d]$. A similar fact holds for $\ell$-adic sheaves (see e.g.\ \cite[Section II.7-8]{KiehlWeissauer}). The above definition mimics this situation insofar as $\SpecZ$ is seen as an analogue of a smooth affine curve. 

Let us give a number of consequences of the preceding axioms, in particular purity, base change and localization: in \refeq{basechange}, suppose that $f$ is smooth and $g: X' \r X$ is a codimension one closed immersion between regular schemes. Then there is a canonical isomorphism 
\eqn \mylabel{eqn_purityExample}
g^! \M_X(Y) = \M_{X'} (X' \x_X Y)(-1)[-2].
\xeqn

Let $Z \subset X$ be a closed immersion of quasiprojective schemes over $S$. Then there is a distinguished triangle of motives with compact support
$$\Mc(Z) \r \Mc(X) \r \Mc(X \backslash Z).$$

Let $S$ be a scheme of equidimension $d$ such that the structural map $f: S \r \SpecZ$ factors as
$$S \stackrel j \lr S' \stackrel i \lr \mathbb A^n_\Z \text{ or } \P[n]_\Z \stackrel p \lr \SpecZ,$$
where $j$ is an open immersion into a regular scheme $S'$, $i$ is a closed immersion and $p$ is the projection map. Then $f^! \one = \one (d-1)[2d-2]$, as one sees by applying relative purity to $p$ and to $j$, and absolute purity to $i$. In particular, the Verdier duality functor on any open subscheme $S$ of $\SpecOF$ is given by $D_{\DMgm(S)} (?) = \IHom(?, \one(1)[2])$ while on $\DMgm(\Fpp)$ it is given by $\IHom(?, \one) = ?\dual$.

\axio \mylabel{axio_Verdierdual} (Verdier dual)
The Verdier dual functor $D$ exchanges ``$!$'' and ``$*$'' throughout, e.g., there are natural isomorphisms $D (f^! M) \cong f^* D(M)$ for any quasi-projective map $f: X \r Y$ and $M \in \DM(Y)$ and similarly with $f_!$ and $f_*$.
\xaxio



\lemm  \mylabel{lemm_reflexivity}
Let $S$ be such that $f^! \one = f^* \one(d)[2d]$ for some integer $d$, where $f: S \r \SpecZ$ is the structural map. For example, $S$ might be regular and affine or projective over $\Z$ (see above), or smooth over $\SpecZ$ (purity). Then, for any compact object $M \in \DMgm(S)$, the canonical map $M \r D(D(M))$ is an isomorphism. This will be referred to as \Def{reflexivity} of Verdier duality. 
\xlemm
\pf
By \refax{compactlygenerated}, it suffices to check it for $M = \pi_* \pi^* \one$, where $\pi: X \r S$ is some map of finite type. In this case it follows for adjointness reasons and the assumption.
\xpf


\axio \mylabel{axio_resolution} (Resolution of singularities) 
Let $K$ be a field. As a triangulated additive tensor category (i.e., closed under triangles, arbitrary direct sums and tensor product), $\DM(K)$ is generated by $\one(-1)$ and all $\M(X)$, where $X / K$ is a smooth projective variety. 

When $S$ is an open subscheme of $\SpecOF$, the generators of $\DM(S)$ are $\one(-1)$, ${i_\pp}_* \M(X_\pp)$, and $\M(X)$, instead, where $X_\pp$ is any projective and smooth variety over $\Fpp$, $i_\pp$ denotes the immersion of any closed point $\Fpp$ of $S$, and $X$ is any regular, flat projective scheme over $\OF$.
\xaxio

Consequently, the subcategories of compact objects $\DMgm(-)$ are generated as a thick tensor  subcategory by the mentioned objects. In Voevodsky's theory of motives over a field of characteristic zero, this is \cite[Section 4.1]{Voevodsky:TCM}. This uses Hironaka's resolution of singularities. Over a field of positive characteristic and number rings, one has to use de Jong's resolution result, see \cite[Lemma B.4]{HuberKahn:Slice}. 

We also need a limit property of the generic point. Let $S$ be an open subscheme of $\SpecOF$, let $\eta: \SpecF \r S$ be the generic point.

\axio \mylabel{axio_generic} (Generic point)
Let $M$ be any geometric motive over $S$. The natural maps $j_* j^* M \r \eta_* \eta^* M$ give rise to an isomorphism $\varinjlim j_* j^* M = \eta_* \eta^* M$, where the colimit is over all open subschemes $j : S' \r S$. It induces a distinguished triangle in $\DM(S)$
\eqn \mylabel{eqn_localgeneric}
\oplus_{\pp \in S} {i_\pp}_* i_\pp^! M \r \id \r \eta_* \eta^* M,
\xeqn
where the sum runs over all closed points $\pp \in S$ and $i_\pp$ is the closed immersion.
\xaxio


\section{Realizations}\mylabel{sect_realizations}

One of the main interests in motives lies in the fact that they are explaining (or are supposed to explain) common phenomena in various cohomology theories. These cohomology functors are commonly referred to as \Def{realization} functors. They typically have the form $\DMgm(S) \r \D^\bound(\C)$, 
where $\C$ is an abelian category whose objects are amenable with the methods of (linear) algebra, such as finite-dimensional vector spaces or finite-dimensional continuous group representations or constructible sheaves. 


For example, let $\ell$ be a prime and let $S$ be either a field of characteristic different from $\ell$ or a scheme of finite type over $\SpecOF[1/\ell]$. The $\ell$-adic cohomology maps any scheme $X$ of finite type over $S$ to 
$$\RG_\ell (X) := \R \pi_* \pi^* \Ql \in \D^\bound_c \left( S, \Ql \right),$$
where $\pi : X \r S$ is the structural map and the right hand category denotes the ``derived'' category of constructible $\Ql$-sheaves on $S$ (committing the standard abuse of notation, see e.g.\ \cite[II.6., II.7.]{KiehlWeissauer}). This functor factors over the \emph{$\ell$-adic realization functor}\index{l-adic realization@$\ell$-adic realization} (\cite[p. 772]{Huber:Realization}, \cite{Ivorra:Realisation}) $\RG_\ell : \DMgm(S) \r \D^\bound_c (S, \Ql)$. When $S$ is of finite type over $\Fpp$, the realization functor actually maps to $\D^\bound_{c, m}(S, \ol \Ql)$, the full subcategory of complexes $C$ in $\D^\bound_c (S, \ol \Ql)$ such that all $\H^n(C)$ are mixed sheaves \cite[1.2]{Deligne:Weil2}.

Further realization functors include Betti, de Rham and Hodge realization. See e.g. \cite[2.3.5]{Huber:Realization}. The following axiom says (in particular) that the $\ell$-adic realization of $\M(X)$ does give the $\ell$-adic cohomology groups.

\axio \mylabel{axio_realfunct} (Functoriality and realizations)
The $\ell$-adic realization functor commutes with the six Grothendieck functors $f_*$, $f_!$, $f^!$, $f^*$, $\t$ and $\IHom$ (where applicable). For example, for any map $f: S' \r S$ and any geometric motive $M$ over $S'$: 
$$(f_* M)_\ell = f_* (M_\ell).$$
\xaxio

\section{Interlude: Perverse sheaves over number rings} \mylabel{sect_perverse}
This section is devoted to a modest extension of $\ell$-adic perverse sheaves \cite{BBD} to the situation where the base $S$ is an open subscheme of $\SpecOF[1/\ell]$. It is needed to formulate \refax{tstructureRealizations} for the $\ell$-adic realization of motives over number rings. This section may be considered a reformulation in ``perverse language'' of the well-known duality for cohomology of the inertia group \cite[Dualit\'e]{SGA4-1/2}. In a nutshell, the theory of perverse sheaves on varieties over $\Fq$ stakes on relative purity, that is $f^! \Zl = f^* \Zl (n)[2n]$ for a smooth map $f$ of relative dimension $n$. The analogous identity for a closed immersion $i: \SpecFpp \r S$ reads
\eqn 
\mylabel{eqn_absolutepurity}
i^! \Zl = i^* \Zl (-1)[-2].
\xeqn 
It is a reformulation of well-known cohomological properties of the inertia group: $\H^1(I_\pp, V) = (V(-1))_{I_\pp}$ for any $\ell$-adic module with continuous $I_\pp$-action ($\pp \nmid \ell$). All higher group cohomologies of $I_\pp$ vanish.    

Let $\D^\bound(S, \Zl)$ be the bounded ``derived'' category of $\Zl$-sheaves on $S$ as constructed by Ekedahl \cite{Ekedahl:Adic}. All following constructions can be done for $\Ql$ instead of $\Zl$, as well. We keep writing $j_*$ for the total derived functor, commonly denoted $\R j_*$ etc. However, $\R^n j_*$ etc.\ keep  their original meaning.

As in \lc, see especially [2.2.10, 2.1.2, 2.1.3, 1.4.10]\footnote{In the sequel, any reference in brackets refers to \cite{BBD}.}, one first defines a notion of stratification, and secondly obtains a $t$-structure on the subcategory $\D^\bound_{(\Sigma, L)}(S, \Zl)$ that are constructible with respect to a given stratification $\Sigma = \{ \Sigma_i \}$ and a set $L$ of irreducible lisse sheaves on the strata. Thirdly, one takes the ``limit'' over the stratifications. The union of all $\D^\bound_{(\Sigma, L)}(S, \Zl)$ is the ``derived'' category $\D^\bound_c(S, \Zl)$ of constructible sheaves. In order to extend the $t$-structure on the subcategories to one on $\D^\bound_c(S, \Zl)$, one has to check that the inclusion $\D^\bound_{(\Sigma', L')}(S, \Zl) \r \D^\bound_{(\Sigma, L)}(S, \Zl)$ is $t$-exact for any refinement of stratifications. Here we employ a different argument. The proof of [2.1.14, 2.2.11] relies on relative purity for $\ell$-adic sheaves \cite[Exp. XVI, 3.7]{SGA4:3}. As in the proof of [2.1.14] we have to check the following: let $\Sigma_i \stackrel{a}\r \Sigma'_i \stackrel b \lr S$ be the inclusions of some strata and let $C \in \D^{\bound, \geq 0}_{(\Sigma', L')} (S, \Zl)$. Then $C \in \D^{\bound, \geq 0}_{(\Sigma, L)} (S, \Zl)$. We can assume $\dim \Sigma_i = 0$, $\dim \Sigma'_i = 1$, since all other cases are clear. Thus, $b$ is an open immersion. We may also assume for notational simplicity that $\Sigma_i = \SpecFpp$. Let $j$ be the complementary open immersion to $a$. By definition, $\H^n b^! C = b^! \H^n C = b^* \H^n C$ is locally constant and vanishes for $n < -1$. In the parlance of Galois modules this means that, viewed as a $\pi_1(\Sigma'_i)$-representation, the action of the inertia group $I_\pp \subset \pi_1(\Sigma'_i)$ on that sheaf is trivial. Thus 
$$a^! \H^n b^* C = a^* (\R^1 j_* j^* \H^n b^* C) [-2] = \H^1 (I_\pp, \H^n b^* C)[-2] = a^* \H^n b^* C (-1)[-2].$$
(We have used $\pp \nmid \ell$ at this point.) The spectral sequence
$$\H^{p-2} a^* \H^q b^! C (-1) = \H^p a^! \H^q b^! C \Rightarrow \H^n a^! b^! C$$
is such that the left hand term vanishes for $p \neq 2$ since $a^*$ is exact w.r.t.~the standard $t$-structure. It also vanishes for $q < -1$ by the above. Hence the right hand term vanishes for $n = p + q < 1$. A fortiori it vanishes for $n < -\dim \Fpp = 0$. 

Objects in the heart of this $t$-structure on $\D^\bound_c(S, \Zl)$ are called \Def{perverse sheaves} on $S$. For example $\Zl[1]$ and $i_* \Zl$ for any  immersion $i$ of a closed point are perverse sheaves on $S$. The \Def{Verdier dual} of any $\C \in \D^\bound_c(S, \Zl)$ is defined by $D(C) := \IHom (C, \Zl(1)[2])$. As above, we have dropped ``$\R$'' from the notation, so that this $\IHom$ means what is usually denoted $\underline {\RHom}$.

\lemm
Let $j: S' \r S$ be an open immersion and $i: Z \r S$ a closed immersion. Let $\eta: \SpecF \r S$ be the generic point. Then $j_*$, $j_!$, $i_*$, $\eta^*[-1]$, $j^*$ and $D$ are $t$-exact, while $i^*$ ($i^!$) is of cohomological amplitude $[-1, 0]$ ($[0, 1]$), in particular right-exact (left- exact, respectively). Finally, the $t$-structure on $\D^\bound_c(S, \Zl)$ is non-degenerate \cite[p. 32]{BBD}. 
\xlemm

\pf
The only non-formal statement is the exactness of $j_*$. The corresponding precursor result [4.1.10] is a reformulation of \cite[Th. 3.1., Exp. XIV]{SGA4:2}, which says for any affine map $j: X \r Y$ over schemes over a field $K$, and any (honest) sheaf $\mathcal F$ which is torsion (prime to $\operatorname{char} K$)
$$d(\R^q j_* \mathcal F)) \leq d(\mathcal F) - q$$
where $d(\mathcal G) := \sup \{ \dim \ol{\{x\}} , \mathcal G_{\ol x} \neq 0 \}$ for any sheaf $\mathcal G$. In our situation, we are given a locally constant sheaf $\mathcal F$ on $S'$ whose torsion is prime to all characteristics of $S$. The conclusion of the theorem also holds for $j$, as follows from the cohomological dimension of $I_\pp$, which is one. 
%
\xpf

Let $\mathcal F$ be any perverse sheaf on $S'$. Following [1.4.22], let the \Def{intermediate extension} $j_{!*} \mathcal F$ be the image of the map $j_! \mathcal F \r j_* \mathcal F$ of perverse sheaves on $S$. As in [2.1.11] one sees that it can be calculated in terms of the good truncation with respect to the standard $t$-structure: $j_{!*} \mathcal F = \tau_{\leq -1}^{can} j_* \mathcal F.$ If $\mathcal F = \mathcal G[1]$, where $\mathcal G$ is a lisse (honest) sheaf on $S'$, this gives $(\R^0 j_* \mathcal G)[1]$.

%

\section{Mixed motives}\mylabel{sect_conjprop}
Throughout this section, let $S = \SpecF$ or $\SpecFpp$ or an open subscheme of $\SpecOF$.\\

This section formulates a number of axioms concerning weights and the motivic $t$-structure on triangulated categories of motives over $S$. In contrast to the axioms listed in \refsect{axiomaticdescription}, the axioms mentioned in this section are wide open, so it might be more appropriate to call them conjectures instead.

\axio \mylabel{axio_cohomdim} (Motivic $t$-structure and cohomological dimension)
The category of geometric motives $\DMgm (S)$ has a non-degenerate $t$-structure \cite[Def. 1.3.1]{BBD} called \Def{motivic $t$-structure}. Its heart is denoted $\MM(S)$\index{MM@$\MM(-)$}. Objects of $\MM(S)$ are called \Def{mixed motives} over $S$. 

For any $M \in \DMgm(S)$, there are $a$, $b \in \Z$ such that $\tau_{\leq a} M = \tau_{\geq b} M = 0$. Here and in the sequel, $\tau_{\leq -}$ and $\tau_{\geq -}$ denote the truncation functors with respect to the motivic $t$-structure.

The \Def{cohomological dimension} of $\DMgm(\Fpp)$ and $\DMgm(F)$ is 0 and 1, respectively, in the sense that
$$\Hom_{\DM(\Fpp)}(M, N[i]) = 0$$
for all mixed motives $M, N$ over $\Fpp$ and $i > 0$ and similarly for mixed motives over $F$ and $i>1$. (For $i < 0$ the term vanishes by the $t$-structure axioms.)

The $t$-structures are such that over $S = \Spec F$ or $\SpecFpp$, $\one \in \MM(S)$, while for an open subscheme $S \subset \SpecOF$, $\one[1] \in \MM(S)$.
\xaxio

The existence of the motivic $t$-structure on $\DMgm(K)$ satisfying the axioms listed in this section is part of the general motivic conjectural framework, see e.g.\ \cite[App. A]{Beilinson:Height}, \cite[Ch. 21]{Andre:Motifs}. The idea of building a triangulated category of motives and descending to mixed motives by means of a $t$-structure is due to Deligne. 
The existence of a motivic $t$-structure on $\DMgm(K)$ is only known in low dimensions: the subcategory of Artin motives, i.e., motives of zero-dimensional varieties, carries such a $t$-structure \cite[Section 3.4.]{Voevodsky:TCM}. By \lc, \cite{Orgogozo:Isomotifs}, the subcategory of $\DMgm(K)$ generated by motives of smooth varieties of dimension $\leq 1$ is equivalent to the bounded derived category of 1-motives \cite[Section 10]{Deligne:Hodge3} up to isogeny. Finally, if $K$ is a field satisfying the Beilinson-Soul\'e vanishing conjecture, such as a finite field or a number field, the category of Artin-Tate motives over $K$ enjoys a motivic $t$-structure \cite{Levine:TateMotives, Wildeshaus:ATM}. The results on Artin-Tate motives are generalized to bases $S$ which are open subschemes of $\SpecOF$ in \cite{Scholbach:Mixed}.

The conjecture about the cohomological dimension is due to Beilinson. A (fairly weak) evidence for this conjecture is the cohomological dimension of Tate motives over $F$ and $\Fpp$, which is one and zero, respectively. This follows from vanishing properties of $K$-theory of these fields. 



The normalization in the last item is merely a matter of bookkeeping, but is motivated by similar shifts in perverse sheaves (\refsect{perverse}). The existence of a motivic $t$-structure is not expected to hold for motives with integral coefficients.

We do not (need to) assume that the canonical functor $\D^\bound(\MM(S)) \r \DMgm(S)$ is an equivalence of categories or, equivalently \cite[Lemma 1.4.]{Beilinson:Derived}, $\Ext^i_{\MM(S)} (A, B) = \Hom_{\DMgm(S)} (A, B[i])$ for all mixed motives $A$ and $B$.

\axio \mylabel{axio_exactness} (Exactness properties) 
Let $S \subset \SpecOF$ be an open subscheme, let $i: \SpecFpp \r \SpecOF$ be a closed point, $j: U \r S$ an open immersion and $\eta: \Spec F \r S$ the generic point. 

Then $j^* = j^!$, $\eta^*[-1]$, $i_*$, $j_*$ and $j_!$ are exact with respect to the motivic $t$-structures on the involved categories of geometric motives. Further, $i^*$ is right-exact, more precisely it maps objects in cohomological degree $0$ to degrees $[-1, 0]$. Dually, $i^!$ has cohomological amplitude $[0, 1]$. Verdier duality $D$ is ``anti-exact'', i.e., maps objects in positive degrees to ones in negative degrees and vice versa.
\xaxio


The axiom is motivated by the same exactness properties in the situation of perverse sheaves over $\SpecOF[1/\ell]$ (\refsect{perverse}). 
The corresponding exactness properties of the above functors on Artin-Tate motives, where the motivic $t$-structure is available, are established in \cite{Scholbach:Mixed}. 

\defi
The cohomology functor with respect to the motivic $t$-structure on $\DMgm(S)$ is denoted $\pH^*$. For any scheme $X / S$, we write
$$\h^i(X, n) := \pH^i \M_S(X)(n).$$
\xdefi

\axio \mylabel{axio_numhom}
Let $X_\eta / F$ be any smooth projective variety. Then numerical equivalence and homological equivalence (with respect to any Weil cohomology) agree on $X_\eta$.
\xaxio

Let either $S$ be a field and let $C$ stand for the $\ell$-adic realization (in case $\text{char } S \neq \ell$), Betti, de Rham or absolute Hodge realization or let $S \subset \SpecOF[1/\ell]$ be an open subscheme and let $C$ be the $\ell$-adic realization. We write  $\RG_C : \DMgm(S) \r \D^\bound(\C)$ for the realization functor, where $\D^\bound(\C)$ is understood as a placeholder of the target category of $C$. (We abuse the notation insofar as that target category is not a derived category in the strict sense when $C$ is the $\ell$-adic realization.) For all realizations over a field, this category is endowed with the usual $t$-structure on the derived category of an exact category, e.g.\ on $\D^\bound_c(K, \Ql)$ for $\ell$-adic realization. When $C$ is the $\ell$-adic realization over an open subscheme $S$ of $\SpecOF[1/\ell]$, we take the perverse $t$-structure on $\D^\bound_c(S, \Ql)$ defined in \refsect{perverse}. Using this, we have the following axiom:

\axio \mylabel{axio_tstructureRealizations} (Exactness of realization functors)
Realization functors $\RG_C$ are exact with respect to the motivic $t$-structure on $\DMgm(S)$. Equivalently, as the $t$-structure on $\D^\bound(\C)$ is non-degenerate, $\RG_C(\pH^0 M) = \pH^0 \RG_C (M)$ for any geometric motive $M$ over $S$. On the left, $\pH^0$ denotes the cohomology functor belonging to the motivic $t$-structure on $\DMgm(S)$, while on the right hand side it is the one belonging to the $t$-structure on $\D^\bound(\C)$.
\xaxio

This axiom is, if fairly loosely, motivated by a similar fact in the theory of mixed Hodge modules: let $X$ be any complex algebraic variety. Then, under the faithful ``forgetful functor'' from the derived category of mixed Hodge modules to the derived category of constructible sheaves with rational coefficients
$$\D^\bound (\category{MHM}(X)) \r \D^\bound_c (X, \Q)$$
the category $\category{MHM}(X)$ corresponds to perverse sheaves on $X$. 


Recall that in an abelian category $\C$, a morphism $f: (X, W^*) \r (Y, W^*)$ between filtered objects is called \Def{strict} if $f(W^n X) = f(X) \cap W^n Y$ for all $n$. 

\axio \mylabel{axio_weight} (Weights)
Any mixed motive $M$ over $S$ has a functorial finite exhaustive separated filtration $W_* M$ called \Def{weight filtration}, i.e., a sequence of subobjects in the abelian category $\MM(S)$
$$0 = W_a M \subset W_{a+1} M \subset \dots \subset W_b M = M.$$ 
Any morphism between mixed motives is strict with respect to the weight filtration.

Tensoring any motive with $\one(n)$ shifts its weights by $-2n$.

Let $\RG_C : \DMgm(S) \r \D^\bound(\C)$ be any realization functor that has a notion of weights (such as the $\ell$-adic realization when $S = \SpecFpp$ or the Hodge realization when $S = \Spec \Q$). Then 
$$\gr_n^W \RG_C (M) = \RG_C (\gr_n^W M)$$
for any mixed motive $M$ over $S$.
\xaxio

\defi
For any $M \in \MM(S)$, we write $\wt(M)$ for the (finite) set of integers $n$ such that $\gr^W_n M \neq 0$. For $M \in \DMgm(S)$,  define\index{wt@$\wt(-)$} $\wt(M) := \cup_{i \in \Z} \wt(\pH^i(M)) -i$.
\xdefi


\axio \mylabel{axio_respectweights} (Preservation of weights)
Let $f: X \r S$ be a quasi-projective map. Then the functors $f_! f^*$ preserve negativity of weights, i.e., given a geometric motive $M$ over $S$ with weights $\leq 0$, $f_! f^* M$ also has weights $\leq 0$. Dually, $f_* f^!$ preserves positive weights. 

In the particular case $S \subset \SpecOF$ (open), let $j: U \r S$ and $\eta : \SpecF \r S$ be an open immersion into $S$ and the generic point of $S$, respectively. Let $i: \SpecFpp \r S$ be a closed point. Then, $i^*$ and $j_!$ preserve negativity of weights and dually for $i^!$ and $j_*$. Finally, $j^*$ and $\eta^*$ both preserve both positivity and negativity of weights.
\xaxio

The preceding weight axioms are motivated by the very same properties of $\ell$-adic perverse sheaves on schemes over $\CC$ or finite fields \cite[5.1.14]{BBD}, number fields \cite{Huber:Perverse} as well as Hodge structures \cite[Th. 8.2.4]{Deligne:Hodge3} and Hodge modules (see \cite[Chapter 14.1]{PetersSteenbrink} for a synopsis). In these settings, actually $f_!$ and $f^*$ preserve negative weights, but we do not need weights for motives over more general bases than the ones above. The weight formalism we require is stronger than the one provided by the differential-graded interpretation of $\DMgm$ over a field \cite{Bondarko:Differential} or \cite[6.7.4]{BeilinsonVologodsky}. 

\bem
Over $S = \SpecOF$, we actually only use the following weight properties: for any $M \in \DMgm(S)$, the interval $\wt(M)$ containing the weights of $M$ satisfies the following two properties: first, it is compatible under functoriality as in \ref{axio_respectweights} and, second, $j_{!*}$ preserves weights of pure smooth motives. (See Definitions \ref{defi_intermediateextension}, \ref{defi_genericallysmooth} for these two notions and the proof of \refth{summaryH1f}.)
\xbem 

\bsp \mylabel{bsp_weightTate}
For any projective (smooth) scheme $X$ of finite type over $S$, the weights of $\h^i(X)(n)$ are $\leq i - 2n$ ($\geq i-2n$, respectively).
\xbsp

\axio \mylabel{axio_mixedpure} (Mixed vs. pure motives)
For any field $K$, the subcategory of pure objects in $\MM(K)$ identifies with $\PureMot_\num(K)$, the category of numerical pure motives over $K$. 
\xaxio

By \refax{cohomdim}, there is an exact sequence 
$$0 \r \H^1(\h^{2n-1}(X_\eta, n)) \r \H^{2n}(X_\eta, n) \r \H^0(\h^{2n}(X_\eta, n)) \r 0.$$
By Axioms \ref{axio_numhom} and \ref{axio_mixedpure}, it reads
\eqn \mylabel{eqn_morphismsChow}
0 \r \CH^n(X_\eta)_{\Q,\hom} \r \CH^n(X_\eta)_\Q \r \CH^n(X_\eta)_\Q / \hom \r 0.
\xeqn
Here $\CH^m(X_\eta)_{\Q,\hom}$ and \index{CHmmodulohom@$\CH^m(-)_\Q / \hom$} $\CH^m (X_\eta)_\Q / \hom$ are by definition the kernel and the image (seen as a quotient of the Chow group) of the cycle class map from the $m$-th Chow group to $\ell$-adic cohomology of $X_\eta$, $\CH^m (X_\eta)_\Q \r \H^{2m}(X_\eta, \Ql(m))$ \cite[VI.9]{Milne:EC}.

As a consequence of the weight filtration, every mixed motive is obtained in finitely many steps by taking extensions of motives in $\PureMot_\num(K)$. Recall also that for any $X / \Fq$ which is smooth and projective the spectral sequence
$$\Ext^p_{\MM(\Fq)}(\one, \h^q(X)) \Rightarrow \Hom_{\DMgm(\Fq)}(\one, \M(X)[p+q])$$
degenerates by \refax{cohomdim} and yields an agreement 
\eqn
\mylabel{eqn_numrat}
\CH^q(X) / \num = \Hom_{\PureMot_\num(\Fq)}(\one, \h_\num^q(X)) \stackrel{\text{\ref{axio_mixedpure}}}= \Hom_{\MM(\Fq)}(\one, \h^q(X)) = \CH^q(X),
\xeqn
i.e., the agreement of rational and numerical equivalence (and thus, of all adequate equivalence relations). 


\bem \mylabel{bem_hardLefschetz}
Recall that the agreement of numerical and homological equivalence on all smooth projective varieties over $F$ implies the motivic hard Lefschetz \cite[5.4.2.1]{Andre:Motifs}: for such a variety $X_\eta / F$ of constant dimension $d_\eta$, let $i \leq d_\eta$ and $a$ any integer. Then, taking the  $(d_\eta - i)$-fold cup product with the cycle class of a hyperplane section with respect to an embedding of $X_\eta$ into some projective space over $F$ yields an isomorphism (``hard Lefschetz isomorphism'')
\eqn \mylabel{eqn_hardLefschetz}
\h^i(X_\eta, a) \stackrel{\cong}{\lr} \h^{2d_\eta - i}(X_\eta, d_\eta - i + a).
\xeqn
The hard Lefschetz is known to imply a non-canonical decomposition \cite{Deligne:Decompositions} 
$$\M(X_\eta) \cong \bigoplus \h^n(X_\eta)[-n].$$
\xbem

We need to assume the following generalization of this. It will be used in \refle{smooth}, which in turn is crucial in \refsect{fcohomology}. Note that the index shift in the second part is due to the normalization in \refax{cohomdim}: for $S = \SpecOF$ and a closed point $i$ as above, take for example $X = S$, $\M(S) = \one = \h^{1}(S)[-1]$ (sic) and $i^* \M(S) = \one_{\Fpp} = \h^0(\Spec \Fpp)$.
  
\axio \mylabel{axio_smoothprojective} (Decomposition of smooth projective varieties) 
Let $X / S$ be smooth and projective. In $\DMgm(S)$, there is a non-canonical isomorphism
$$\phi_X : \M(X) \cong \bigoplus_n \h^n(X)[-n].$$
For open subschemes $S \subset \SpecOF$, this isomorphism is compatible with pullbacks along all closed points $i: \SpecFpp \r S$ in the following sense: let $X_\pp$ be the fiber of $X$ over $\Fpp$, and let $\psi$ be the isomorphism making the following diagram commutative. Its left hand isomorphism is an instance of base change.
$$\xymatrix{
i^* \M(X) \ar[r]^(.4){i^* \phi_X} \ar[d]^\cong & 
\oplus_n i^* \h^n(X)[-n] \ar[d]^\psi \\
\M(X_\pp) \ar[r]^(.4){\phi_{X_\pp}} &
\oplus_m \h^m(X_\pp)[-m]}
$$
Then $\psi$ respects the direct summands, i.e., induces isomorphisms
$$i^* \h^n(X)[-n] \cong  \h^{n-1}(X_\pp)[-n+1].$$
\xaxio

\section{Motives over number rings}\mylabel{sect_motivesOF}

In the following sections we assume the axioms of Sections \ref{sect_axiomaticdescription}, \ref{sect_realizations}, and \ref{sect_conjprop}. Unless explicitly mentioned otherwise, let $S$ be an open subscheme of $\SpecOF$, let $i: \SpecFpp \r \SpecOF$ be a closed point, $j: S' \r S$ an open subscheme and $\eta: \Spec F \r S$ the generic point.

This section derives a number of basic results about motives over $S$ from the axioms spelled out above. We define and study the \emph{intermediate extension} $j_{!*}: \MM(S') \r \MM(S)$ in analogy to perverse sheaves (\refde{intermediateextension}). An ``explicit'' set of generators of $\DMgm(S)$ (\refsa{generatorsoverOF}) is obtained using $j_{!*}$. We introduce a notion of \emph{smooth motives} (\refde{genericallysmooth}), which should be thought of as analogs of lisse sheaves. Using this notion, we extend the intermediate extension to a functor $\eta_{!*}$ spreading out motives over $F$ with a certain smoothness property to motives over $S$, cf.\  \refde{genericintermediateextension}. This functor will be the main technical tool in dealing with $f$-cohomology in \refsect{fcohomology}. In \refle{excreal} we express the $\ell$-adic realization of motives of the form $j_{!*} M$ in sheaf-theoretic language.

\subsection{Cohomological dimension}

The following is an immediate consequence of \refax{exactness}:
\lemm \mylabel{lemm_eta}
For any scheme $X$ over $S$ we have $\eta^* [-1] \h^i(X,n) = \h^{i-1}(X \x_S F, n)$. 
\xlemm

The following lemma parallels (and is a consequence of) \refax{cohomdim}. 

\lemm \mylabel{lemm_cohomdimOF} The cohomological dimension of $\DMgm(S)$ is two, that is to say, for any two mixed motives $M, N$ over $S$,
$$\Hom_{\DMgm(S)}(M, N[i]) = 0$$
for all $i > 2$. In particular $\H^{i}(M)$ vanishes for $|i|>1$. 
\xlemm
\pf
Apply $\Hom(M,-)$ to the localization triangle $\oplus_{\pp \in S} {i_\pp}_* i_\pp^! N \r N \r \eta_* \eta^* N$ of \refax{generic}, where $i_\pp$ are the immersions of the closed points of $S$. The terms adjacent to $\Hom(M, N[i])$ are $\Hom(M, \oplus_{\pp} {i_\pp}_* i_\pp^! N[i]) = \oplus_{\pp} \Hom(i_\pp^* M, i_\pp^! N[i])$ (as $M$ is compact) and $\Hom(M, \eta_* \eta^* N[i]) = \Hom(\eta^* M, \eta^* N[i])$. The latter term vanishes for $i>1$ since $\eta^*[-1]$ is exact and the cohomological dimension of $\DMgm(F)$ is one.

To deal with the former term, we have to take into account that $i_\pp^!$ and $i_\pp^*$ are not $t$-exact, but of cohomological amplitude $[0, 1]$ and $[-1, 0]$, respectively. By decomposing $i_\pp^! N$ into its $\pH^{1}$- and $\pH^0$-part and similarly with $i_\pp^* M$ and using that the cohomological dimension of $\DMgm(\Fpp)$ is zero, the term vanishes for $i > 2$. Using general $t$-structure properties, the second claim is a particular case of the first one.
\xpf

\subsection{Intermediate extension}
\defi (Motivic analog of \cite[Def. 1.4.22]{BBD}) \mylabel{defi_intermediateextension}
The \Def{intermediate extension} $j_{!*}$\index{j@$j_{"!*}$} of some mixed motive $M$ over $S'$ is defined as
$$j_{!*} M := \Im (j_! M \r j_* M).$$
The image is taken in the abelian category $\MM(S)$, using the exactness of $j_!$ and $j_*$, \refax{exactness}.
\xdefi

\bem \mylabel{bem_triangles}
Let $i: Z \r S$ be the complement of $j$.  The localization triangles (\refax{localization}) and cohomological amplitude of $i^*$ (\refax{exactness}) yield short exact sequences in $\MM(S)$
\eqn \mylabel{eqn_local1}
0 \r i_* \pH^{-1} i^* j_* M \r j_! M \r j_{!*} M \r 0,
\xeqn
\eqn \mylabel{eqn_local2}
0 \r j_{!*} M \r j_* M \r i_* \pH^{0} i^* j_* M \r 0.
\xeqn
These triangles are the same as for perverse sheaves in the situation that the analog of \refax{exactness}, \cite[4.1.10]{BBD}, is applicable.
\xbem

\lemm \mylabel{lemm_properties}
\begin{itemize}
\item 
Given any mixed motive $M$ over $S'$, $j_{!*} M$ is, up to a unique isomorphism, the unique mixed extension of $M$ (i.e., an object $X$ in $\MM(S)$ such that $j^* X = M$) not having nonzero subobjects or quotients of the form $i_* N$, where $i: Z \r S$ is the closed complement of $j$ and $N$ is a mixed motive on $Z$. 
\item
For any two composable open immersions $j_1$ and $j_2$ we have ${j_1}_{!*} \circ {j_2}_{!*} = ({j_1} \circ j_2)_{!*}$.
\item
$j_{!*}$ commutes with duals, i.e., $D(j_{!*} -) \cong j_{!*} D(-)$.
\end{itemize}
\xlemm
\pf
The proofs of the same facts for perverse sheaves \cite[Cor. 1.4.25, 2.1.7.1]{BBD} carry over to this setting. 
The first statement easily implies the last one.
\xpf

The following proposition makes precise the intuition that any motive $M$ over $S$ should be reconstructed by its generic fiber (over $F$) and a finite number of special fibers (over various $\Fpp$).

\satz \mylabel{satz_generatorsoverOF}
As a thick subcategory of $\DM(S)$, $\DMgm(S)$ is generated by motives of the form
\begin{itemize}
\item 
$i_* \M(X_\pp)(m)$, where $X_\pp / \Fpp$ is smooth projective, $m \in \Z$ and $i: \SpecFpp \r S$ is any closed point and 
\item 
$j_{!*} j^* \h^k(X, m)$, where $X$ is regular, flat and projective over $S$ with smooth generic fiber, and $j: S' \r S$ is such that $X \x_S S'$ is smooth over $S'$ and $k$ and $m$ are arbitrary. 
\end{itemize}
\xsatz
\pf 
Let $\mathcal D \subset \DMgm(S)$ be the thick category generated by the objects in the statement. 
By resolution of singularities over $S$ (\refax{resolution}), $\DMgm(S)$ is the thick subcategory of $\DM(S)$ generated by objects $i_* \M(X_\pp)(m)$ and $\M(X)(m)$, where $X_\pp$ and $X$ are as in the statement and $m \in \Z$. 

It is therefore sufficient to see $M := \M(X) \in \mathcal D$. Let $j : S' \r S$ be such that $X_{S'}$ is smooth over $S'$. By \ref{axio_localization} it is enough to show $j_* j^* M \in \mathcal D$, since motives over finite fields are already covered. Applying the truncations with respect to the motivic $t$-structure to $j_* j^* M$ and exactness of $j_*$, $j^*$ (\refax{exactness}) shows that we may deal with $j_* j^* \h^k(X, m)$ for all $k$ instead of $j_* j^* M$. (Only finitely many $k$ yield a nonzero term by \refax{cohomdim}.) By \refbe{triangles}, there is a short exact sequence of mixed motives 
$$0 \r j_{!*} j^* \h^k(X, m) \r j_* j^* \h^k(X, m) \r i_* \pH^0 i^* j_* j^* \h^k(X, m) \r 0.$$
Here $i$ is the complement of $j$. The left and right hand motives are in $\mathcal D$, hence so is the middle one.
\xpf

\subsection{Smooth motives} \mylabel{sect_smoothmotives}
The notion of smooth motives (a neologism leaning on lisse sheaves) is a technical stepstone for the definition of the generic intermediate extension $\eta_{!*}$, cf.\ \refde{genericintermediateextension}. Roughly speaking, smoothness for mixed motives $M$ means that $i^* M $ and $i^! M$ do not intermingle in the sense that their cohomological degrees are disjoint.

\defi \mylabel{defi_genericallysmooth}
Let $M$ be a geometric motive over $S$. It is called \Def{smooth}\index{motive!smooth} if for any closed point $i: \SpecFpp \r S$ there is an isomorphism
$$i^! M \cong i^* M (-1)[-2].$$
$M$ is called \Def{generically smooth}\index{motive!generically smooth} if there is an open (non-empty) immersion $j: S' \r S$ such that $j^* M$ is smooth.
\xdefi

Let $X / S$ be a scheme with smooth generic fiber $X_\eta$. Then $\M_S(X)$ is a generically smooth motive. 

The isomorphism in \refde{genericallysmooth} is not required to be canonical in any sense. Therefore, the subcategory of smooth motives is \emph{not} triangulated in $\DMgm(S)$. 

\lemm \mylabel{lemm_basic}
Let $M$ be a smooth mixed motive over $S$. Let $i: Z \r S$ be proper closed subscheme, let $j: S' \r S$ be its complement. Then $i^! M = (\pH^1 i^! M)[-1]$ and dually $i^* M = (\pH^{-1} i^* M)[1]$.
\xlemm
\pf
By assumption $i^! M \cong i^* M(-1)[-2]$. By \refax{exactness}, the left hand side of that isomorphism is concentrated in degrees $[0, 1]$. The right hand side is in degrees $[1, 2]$. This shows $i^! M = \pH^1 (i^! M)[-1]$ by \refax{cohomdim} and similarly for $i^* M$.
\xpf

The following is the key relation of smooth motives and the intermediate extension. Note the similarity with \refle{locconst}.

\lemm \mylabel{lemm_smooth_and_j}
Let $M$ be a smooth mixed motive over $S$. Then $M$ is canonically isomorphic to $j_{!*} j^* M$.
\xlemm
\pf
Let $i : Z \r S$ be the complement of $j$. Given any $i_* N \subset M$ with $N \in \MM(Z)$, we apply the left-exact functor $i^!$ and see $N \subset \pH^0(i^! M) \stackrel{\ref{lemm_basic}}= 0$. Quotients of $M$ of the form $i_* N$ are treated dually. We now invoke \refle{properties}.
\xpf

\lemm \mylabel{lemm_smooth}
Let $X$ be any smooth projective scheme over $S$. Set $M := \M(X)$. Then all $\h^n X = \pH^n M$ are smooth.
\xlemm
\pf
Let $f_{m,n}$ be the $(m,n)$-component of the bottom isomorphism making the following commutative:
$$
\xymatrix{
i^! M \ar[r]^{\cong, \text{ see \refeq{purityExample}}} \ar[d]^{\cong, \ref{axio_smoothprojective}} & 
i^* M(-1)[-2] \ar[d]^{\cong, \ref{axio_smoothprojective}} \\
\oplus_m A_m := \oplus i^! (\pH^m M)[-m] \ar[r]^{\cong} & 
\oplus_n B_n := \oplus i^* (\pH^n M)(-1)[-n-2].
}$$
We claim $f_{m,n} = 0$ for all $m \neq n$, from which the lemma follows. By \refax{smoothprojective} we have $B_n \cong \h^{n-1}(X_\pp)[-n-1](-1)$. Using this and the reflexivity of the Verdier dual functor, we obtain an isomorphism $A_m \cong (\pH^{m+1} i^! M)[-1-m]$. Hence $B_n$ is concentrated in cohomological degree $n+1$, while $A_m$ is in degree $n+2$. (The a priori bounds of \refax{exactness} would be $[m, m+1]$ and $[n+1, n+2]$, respectively.) As the cohomological dimension of motives over $\Fpp$ is zero (\refax{cohomdim}), the only way for $f_{m,n} \neq 0$ is $m=n$. 
\xpf

\subsection{Generic intermediate extension}
\lemm \mylabel{lemm_spreadout} (Spreading out morphisms)
Given two geometric motives $M$ and $M'$ over $S$ together with a map $\phi_\eta: \eta^* M \r \eta^* M'$, there is an open subscheme $j: S' \subset S$ and a map $\phi_{S'}: j^* M \r j^* M'$ which extends $\phi_\eta$. Any two such extensions agree when restricted to a possibly smaller open subscheme. In particular, if $\phi_\eta$ is an isomorphism, then $\phi_{S'}$ is an isomorphism for sufficiently small $S'$.
\xlemm
\pf
The adjunction map $M \r \eta_* \eta^* M$ and $\eta_* \phi_\eta$ give a map $M \r \eta_* \eta^* M'$, hence by \refeq{localgeneric} a map $M \r \oplus_\pp {i_\pp}_* i_\pp^! M'[1]$. As $M$ is compact, it factors over a finite sum $\oplus_{\pp \in T} {i_\pp}_* i_\pp^! M'[1]$. Let $j: S' \r S$ be the complement of the points in $T$. The map $M \r \eta_* \eta^*M'$ factors over ${j}_* j^* M'$ and gives a map $j^* M \r j^* M'$ which extends $\phi_\eta$. The first claim is shown. 

For the unicity of the extension, we may assume that $\phi_\eta$ is zero, and show that $\phi_{S'}$ is zero for some suitable $S'$. This is the same argument as before: the map $M \r {j}_* j^* M'$ constructed in the previous step factors over $\oplus_{\pp \in S'} {i_\pp}_* i_\pp^! M'$, since $M \r \eta_* \eta^* M'$ is zero. By compacity of $M$, only finitely many primes in the sum contribute to the map, discarding these yields the claim.

If $\phi_\eta$ is an isomorphism, $\psi_\eta := \phi_\eta^{-1}$ can be extended to some $\psi_{S'}$. As both $\phi_{S'} \circ \psi_{S'}$ and $\id_{S'}$ extend $\id_\eta$, they agree on some possibly smaller open subscheme of $S$ and similarly with $\psi_{S'} \circ \phi_{S'}$.
\xpf

\bem
The lemma shows the full faithfulness of the functor 
$$\varinjlim_{S' \subset S} \DMgm(S') \stackrel{\eta^*}{\lr} \DMgm(F).$$
Its essential surjectivity is a consequence of \refax{compactlygenerated}, so we have an equivalence. However, we will stick to the more basic language of colimits in $\DM(S)$ instead of colimits of the categories of geometric motives.
\xbem

\defi \mylabel{defi_genericintermediateextension}
Let $M_\eta \in \DMgm(F)$ be a motive such that there exists a generically smooth mixed motive $M$ over $S$ (\refde{genericallysmooth}) with $\eta^* M \cong M_\eta$. Then the \Def{generic intermediate extension} $\eta_{!*} M_\eta$ is defined as 
$$\eta_{!*} M_\eta := j_{!*} j^* M$$
where $j: S' \r S$ is an open immersion such that $j^* M$ is smooth. 
\xdefi

This is independent of the choices of $j$ and $M$ (Lemmas \ref{lemm_smooth_and_j}, \ref{lemm_spreadout}) and functorial (\ref{lemm_spreadout}). For a mixed, non-smooth motive $M$, there need not be a map $j_{!*} j^* M \r M$. Therefore, $\varinjlim j_{!*} j^* M$ does not make sense unless there is some smoothness constraint on $M_\eta$.




\subsection{Intermediate extension and $\ell$-adic realization} \mylabel{sect_intermediateell}
This subsection deals with the interplay of the (generic) intermediate extension functor on mixed motives and the $\ell$-adic realization. In this subsection, $S$ is an open subscheme of $\SpecOF[1/\ell]$. The following lemma is well-known.

\lemm \mylabel{lemm_locconst}
Let $\mathcal F$ be an \'etale (honest) locally constant sheaf on $S$. Let $\eta: \SpecF \r \SpecOF[1/\ell]$ be the generic point. Then the canonical map $\mathcal F \r \R^0 \eta_* \eta^* \mathcal F$ is an isomorphism. 
\xlemm

\lemm \mylabel{lemm_excreal}
Let $M$ be a mixed motive over $S'$. Let $j: S' \r S$ be an open immersion. Then
$$(j_{!*} M)_\ell = j_{!*} (M_\ell).$$
Let $i$ be the complementary closed immersion to $j: S' \r S$ and let $\eta'$ and $\eta$ be the generic point of $S'$ and $S$, respectively. If $M$ is additionally smooth, one has
$$(i^* j_{!*} M)_\ell = i^* j_{!*} M_\ell = i^* (\R^0 \eta_* \eta'^* M_\ell[-1])[1].$$
\xlemm

To clarify the statement, note that the $\ell$-adic realization of $M$ is a perverse sheaf on $S'$ by \refax{tstructureRealizations}. Thus, $j_{!*}$ (\refsect{perverse}) can be applied to it.

\pf
The first statement follows from \refax{realfunct}, the definition of $j_{!*}$ and the exactness of $\RG_\ell$ (\refax{tstructureRealizations}).
%

Let now $M$ be mixed and smooth over $S'$. As $M_\ell$ is a perverse sheaf by \ref{axio_tstructureRealizations}, there is an open immersion $j' : S'' \r S'$ such that $j'^* M_\ell[-1]$ is a locally constant (honest) sheaf on $S''$. As $M$ is smooth we know from Lemmas \ref{lemm_properties} and \ref{lemm_smooth_and_j}
$$i^* j_{!*} M = i^* (j \circ j')_{!*} j'^* M.$$
By the interpretation of $(j \circ j')_{!*}$ in terms of $\R^0 (j \circ j')_*$ (\refsect{perverse}) we have 
$$(i^* j_{!*} M)_\ell = i^* j_{!*} M_\ell = i^* (\R^0  (j \circ j')_* j'^* M_\ell[-1])[1] \stackrel{\text{\ref{lemm_locconst}}}= i^* (\R^0 \eta_* \eta'^* M_\ell[-1])[1].$$
\xpf

\section{$f$-cohomology} \mylabel{sect_fcohomology}

\subsection{$f$-cohomology via non-ramification}
Let $F$ be a number field. For any place $\pp$ of $F$, let $F_\pp$ be the completion, $G_\pp$ the local Galois group. For finite places, $I_\pp$ denotes the inertia group. For brevity, we will usually write $\H^*(M)$ for $\H^*(S, M)$, where $M$ is any motive over some base $S$.

\defi \cite[Section 3]{BlochKato} \mylabel{defi_H1flocal}
Let $V$ be a finite-dimensional $\ell$-adic vector space, endowed with a continuous action of $G_\pp$, where $\pp$ is a finite place of $F$ not over $\ell$. Set
$$\H^i_f(F_\pp, V) := \left\{ \begin{array}{ll}
     \H^0 (F_\pp, V) & i=0\\
\ker \H^1(F_\pp, V) \r \H^1(I_\pp, V) & i=1 \\
0 & \text{else.}
   \end{array} \right.$$  
\xdefi

\bem
If $\pp$ lies over $\ell$, the definition is completed by $\H^1_f(F_\pp, V) := \ker \H^1(F_\pp, V) \r \H^1(F_\pp, B_{crys} \t V)$, where $B_{crys}$ denotes the ring of $\pp$-adic periods \cite{FontaineMessing}. We will disregard this case throughout.
\xbem

\lemm \mylabel{lemm_charH1f}
Let $\eta_\pp: \Spec F_\pp \r \Spec \OFpp$ be the generic point of the completion of $\OF$ at $\pp$. Using the above notation, for $\pp$ not over $\ell$, there is a canonical isomorphism $\H^1_f (F_\pp, V) \cong \H^1 (\OFpp, \R^0 {\eta_\pp}_{*} V)$. (The right hand side denotes $\ell$-adic cohomology over $\OFpp$.)
\xlemm
\pf
For any $\ell^n$-torsion sheaf $\mathcal F$ on $F_\pp$ we write $A(\mathcal F) := \ker \H^1(F_\pp, \mathcal F) \r \H^1(I_\pp, \mathcal F)$. The $\Ql$-sheaf $V$ is, by definition, of the form $U \t_\Zl \Ql$, where $U = (U_n)_n$ is a projective system of $\Z / \ell^n$-sheaves. By definition
$$\H^1(F_\pp, V) = \varprojlim_{n \in \N} \H^1(F_\pp, U_n) \t \Ql$$
and similarly for $\H^1(I_\pp, V)$. Both $\varprojlim_n$ and $- \t_\Zl \Ql$ are left-exact functors, so one has
$$\H^1_f(F_\pp, V) = \left (\varprojlim_n A(U_n) \right) \t \Ql.$$
Thus it is sufficient to show $A(U) = \H^1(\OFpp, \R^0 {\eta_\pp}_{*} U)$ for any $\ell^n$-torsion sheaf $U$ over $F_\pp$.


Recall the description of \'etale sheaves on $\OFpp$ from \cite[II.3.12, II.3.16]{Milne:EC}. Let $i: \SpecFpp \r \Spec \OFpp$ be the closed point. As $\OFpp$ is a henselian ring \cite[Prop. I.4.5]{Milne:EC}, for any sheaf $\mathcal F$ on $\Spec \OFpp$, the global sections depend only on the special fiber and
$$ \Gamma_{\Spec F_\pp} =  \Gamma_{\Spec \OFpp} \circ ({\eta_\pp}_{*}) = \Gamma_{\Spec \OFpp} \circ (i_* i^* {\eta_\pp}_{*}). $$
These functors can be interpreted using group cohomology: $\Gamma_{\Spec \OFpp} \circ i_* = \Gamma_{\Fpp}$ and $(-)^{I_\pp} = i^* {\eta_\pp}_{*}$ (\lc). The Hochschild-Serre spectral sequence for $(-)^{G_\pp} = (-)^{\Gal(\Fpp) } \circ (-)^{I_\pp}$ can be translated to 
$$ \H^p ( \Spec \OFpp, i_* i^* \R^q {\eta_\pp}_* U) \Rightarrow \H^n (F_\pp, U).$$
In addition we have the Leray spectral sequence
$$ \H^p ( \Spec \OFpp, \R^q {\eta_\pp}_* U) \Rightarrow \H^n (F_\pp, U). $$
The exact sequence of low degrees of the Hochschild-Serre sequence maps to the sequence below:
$$
\xymatrix
{
0 \ar[r] &
\H^1(\Spec \OFpp, \R^0 {\eta_\pp}_{*} U) \ar[r] &
\H^1(F_\pp, U) \ar[r] \ar[d]^= &
\H^0(\Spec \OFpp, \R^1 {\eta_\pp}_* U) \ar[d]
\\
0 \ar[r] &
A(U) \ar[r] &
\H^1(F_\pp, U) \ar[r] &
\H^1(I_\pp, U)
}
$$
As $\H^0(\Gal(\Fpp), \H^1(I_\pp, U)) \subset \H^1(I_\pp, U)$ and $\Gamma_\OFpp = \Gamma_\OFpp \circ i_* i^*$, the right hand map is injective, therefore there is a unique isomorphism between the left hand terms making the diagram commutative.
\xpf

In order to proceed to a global level, the following definition is done:

\defi \mylabel{defi_H1fglobal} \cite[II.1.3]{FPR}
Given an $\ell$-adic continuous representation $V$ of $G = \Gal(F)$, define $\H^i_f(F, V)$\index{H1fF@$\H^1_f(F, -)$} to be such that the following diagram is cartesian. In the lower row, $V$ is considered a $G_\pp = \Gal(F_\pp)$-module by restriction. 
\begin{displaymath}
\xymatrix
{
\H^i_f(F, V) \ar[r] \ar[d] &
\H^i(F, V) \ar[d] \\
\prod \H^i_f(F_\pp, V) \ar[r] &
\prod \H^i(F_\pp, V)
}
\end{displaymath}

The product ranges over all finite places $\pp$ of $F$. 
We define $\H^i_{f,\backslash \crys}(F, V)$\index{H1flF@$\H^i_{f,\backslash \ell}(F, V)$} similarly, except that in the lower row of the preceding diagram only places $\pp$ that do not lie over $\ell$ occur.
\xdefi

\lemm \mylabel{lemm_charH1fglobal}
Let $V$ be an $\ell$-adic \'etale sheaf on $\Spec F$. Then there is a natural isomorphism
$$\H^i_{f,\backslash \crys}(F, V) \cong \H^1 (\OF[1/\ell], \R^0 {\eta}_{*} V).$$
\xlemm
\pf
By the same argument as in the previous proof, we may assume that $V$ is a sheaf of $\Z / \ell^n$-modules, since the isomorphism we are going to establish is natural in $V$ and 
$$\H^i_{f,\backslash \crys}(F, V) = \ker \H^i(F, V) \r \prod_{\pp \nmid \ell} \left ( \H^i(F_\pp, V) / \H^i_f(F_\pp, V) \right).$$
Consider the following cartesian diagram ($\pp \nmid \ell$)
$$\xymatrix{
\Spec F_\pp \ar[r]^{\eta_\pp} \ar[d]^b &
\Spec \OFpp \ar[d]^a &
\Spec \Fpp \ar[l]^{i_\pp} \ar[d]^= \\
\Spec F \ar[r]^{\eta} &
\Spec \OF[1/\ell] &
\Spec \Fpp \ar[l]^{i} 
}
$$
In the derived category of $\Z / \ell^n$-sheaves on $\SpecOF[1/\ell]$, there is a triangle $\R^0 \eta_* V \r \R \eta_* V \r \R^1 \eta_*[-1] V$. Likewise,
$\R^0 {\eta_\pp}_* b^* V \r {\R \eta_\pp}_* b^* V \r \R^1 {\eta_\pp}_* b^* V[-1]$. (We have used $\pp \nmid \ell$, since the inertia group has cohomological dimension bigger than one for $\pp | \ell$.)
This yields exact horizontal sequences, the vertical maps are adjunction maps

\scriptsize
$$
\xymatrix
{
0 \ar[r] &
\H^1(\Spec \OF[1/\ell], {\eta}_{*} V) \ar[r] \ar[d] &
\H^1(F, V) \ar[r] \ar[d] &
\H^0(\Spec \OF[1/\ell], \R^1 \eta_* V) \ar[d]^\alpha
\\
0 \ar[r] &
\prod_{\pp \nmid \ell} \H^1(\OFpp, \R^0 {\eta_\pp}_* b^* V) \ar[r] &
\prod_{\pp \nmid \ell} \H^1(F_\pp, b^* V) \ar[r] &
\prod_{\pp \nmid \ell} \H^0(\OFpp, \R^1 {\eta_\pp}_* b^* V)
}
$$
\normalsize
We will show that $\alpha$ is injective. Hence, the left square is cartesian and by definition and \refle{charH1f} the claim is shown. Indeed, $\alpha$ factors as 
$$\H^0(\OF[1/\ell], \R^1 \eta_* V) \subset \prod_{\pp \nmid \ell} \H^0(\Fpp, i_\pp^* \R^1 \eta_* V) \r \prod_{\pp \nmid \ell} \H^0(\OFpp, \R^1 {\eta_\pp}_* b^* V) \left ( \stackrel \cong = \prod_{\pp \nmid \ell} \H^0 (\Fpp, i_\pp^* \R^1 {\eta_\pp}_* b^* V) \right).$$
using $i^* \R^1 \eta_* V = i_\pp^* a^* \R^1 \eta_* V = i_\pp^* \R^1 {\eta_\pp}_* b^* V$.
\xpf

\defi \cite[Remark 4.0.1.b]{Beilinson:Height}, \cite[Conj. 5.3]{BlochKato}, \cite[Section 6.5]{Fontaine:Valeurs}, \cite[III.3.1.3]{FPR} \mylabel{defi_H1fEllAdic}
Let $M_\eta$ be a mixed motive over $F$. Let, similarly to \refde{H1fglobal}, $\H^i_f (M_\eta)$ be defined such that the following diagram, in which the bottom products are taken over all primes $\ell$, is cartesian. As usual, ${M_\eta}_\ell$ is the $\ell$-adic realization, seen as a $G$-module.

\begin{displaymath}
\xymatrix
{
\H^i_f(F, {M_\eta}) \ar[r] \ar[d] &
\H^i(F, {M_\eta}) \ar[d] \\
\prod_\ell \H^i_f(F, {M_\eta}_\ell) \ar[r] &
\prod_\ell \H^i(F, {M_\eta}_\ell)
}
\end{displaymath}
Again, to rid ourselves from crystalline questions at $\pp | \ell$, we define $\H^i_{f,\backslash \crys}(F, {M_\eta})$ by replacing $\prod_\ell \H^i_f(F, {M_\eta}_\ell)$ in the bottom row by $\prod_\ell \H^i_{f,\backslash \crys}(F, {M_\eta}_\ell)$.
\xdefi

We are now going to exhibit an interpretation of $f$-cohomology thus defined in terms of the generic intermediate extension $\eta_{!*}$. Recall that we are assuming in this section the axioms of Sections \ref{sect_axiomaticdescription}, \ref{sect_realizations}, and \ref{sect_conjprop}. Mixed motives are needed to even define $\eta_{!*}$. Moreover, for the comparison result, we need to assume the following conjecture. 
 


\lemm \mylabel{lemm_injectivity}
Let $N$ be any mixed motive over $\Fpp$. 
The $\ell$-adic realization map $\H^0(\Fpp, N) \r \H^0(\Fpp, N_\ell) := N_\ell^{\Gal(\Fpp)}$ is injective.
\xlemm
\pf
By the strictness of the weight filtration, the canonical maps 
$$\H^0 (\gr_0^W N) \leftarrow \H^0(W_0 N) \r \H^0(N)$$ 
are both isomorphisms. Moreover, the $\ell$-adic realization functor commutes with $\gr_0^W$ by \refax{weight}, so that we can replace $N$ by $\gr_0^W$ and assume that $N$ is pure of weight $0$. In view of our assumptions on motives, cf.\ \refeq{numrat}, all adequate equivalence relations agree, so that we may regard $N$ as a pure motive with respect to any adequate equivalence relation. As the injectivity is stable under taking direct summands, we may assume $N = h(X, n)$ for $X$ smooth and projective over $\Fpp$, by definition of pure motives and \refax{mixedpure}. The left hand side is given by $\CH^n(X)$, so the map is injective by \refeq{numrat}.
\xpf

\theo \mylabel{theo_charH1fmotivic} 
Let $M$ be a generically smooth mixed motive over $\OF$ (\refde{genericallysmooth}). Set $\eta^* M [-1] =: M_\eta$. 
There is a natural isomorphism
$$\H^0(\OF, \eta_{!*} \eta^* M) \stackrel{\cong}{\lr} \H^1_{f, \backslash \crys}(F, M_\eta).$$
\xtheo
\pf
Notice that $\eta_{!*} \eta^* M$ is well-defined by the assumptions. We want to show that there is a cartesian commutative diagram
$$
\xymatrix
{
\H^0(\eta_{!*} \eta^* M) \ar[r] \ar@{.>}[d] &
\H^0(\eta_* \eta^* M)  = \H^1(M_\eta) \ar[d] 
\\
\prod_\ell \H^1_{f,\backslash \ell}(F, {M_\eta}_\ell) \ar[r] &
\prod_\ell \H^1(F, {M_\eta}_\ell)
}
$$
Let $j: U \r \SpecOF$ be any open immersion  such that $j^* M$ is smooth. We have $\eta_{!*} \eta^* M = j_{!*} j^* M$. The left hand term of the exact sequence 
$$\oplus_{\pp \in U} \H^0({i_\pp}_* i_\pp^! M) \r \H^0(j_* j^* M) \r \H^0(\eta_* \eta^* M) \r \oplus_{\pp \in U} \H^1({i_\pp}_* i_\pp^! M)$$
induced by \refeq{localgeneric} vanishes as $i_\pp^! M$ is concentrated in cohomological degree $1$ for $\pp \in U$ (\refle{basic}). Any $a \in \H^0(\eta^* M)$ maps to a finite sub-sum of $\oplus_{\pp \in \SpecOF} \H^1({i_\pp}_* i_\pp^! M)$, so letting $j$ be the open complement of these points, $a$ lies in (the image of) $\H^0(j_* j^* M)$:
$$\H^0(\eta^* M) = \varinjlim_{j: U \r \SpecOF \atop j^* M \text{ smooth}} \H^0(j_* j^* M).$$
By \refle{prepinjectivity} below, the map $\H^0 (j_{!*} j^* M) \r \H^0(j_* j^* M) \r \H^0(\eta^*M)$ is injective. Therefore, taking the colimit over all $U$ such that $M|_U$ is smooth, the exact localization sequence 
$$0 \r \H^0(j_{!*} j^* M) \r \H^0(j_* j^* M) \r \oplus_{\pp \notin U} \H^0 (\pH^0 i_\pp^* j_* j^* M)$$
stemming from \refeq{local2} gives
$$0 \r \H^0(j_{!*} j^* M) \r \H^0(\eta_* \eta^* M) \r \oplus_{\pp} \H^0 (\pH^0 i_\pp^* {j_\pp}_* j_\pp^* M).$$
Here $j_\pp$ is the complementary open immersion to $i_\pp$ and the direct sum is over all (finite) places $\pp$ of $\OF$. We have $i_\pp^* \eta_* \eta^* M = i_\pp^* {j_\pp}_* j_\pp^* M$, so the top sequence in the following commutative diagram is exact:
\scriptsize
\eqn \mylabel{eqn_diagramreal}
\xymatrix{
0 \ar[r] & 
\H^0(j_{!*} j^* M) \ar[r] \ar[d] & 
\H^1(M_\eta) \ar[r] \ar[d] &
 \oplus_{\pp} \H^0(\pH^0 i_\pp^* \eta_* \eta^* M) \ar[d] \\
0 \ar[r] &
\prod_\ell \H^0((j_\ell^* j_{!*} j^* M)_\ell) \ar[r] &
\prod_\ell \H^1({M_\eta}_\ell) \ar[r] &
\prod_\ell \oplus_{\pp \nmid \ell} \H^0((\pH^0 i_\pp^* \eta_* \eta^* M)_\ell)
}
\xeqn
\normalsize
The lower row  denotes $\ell$-adic cohomology over $\OF[1/\ell]$, $F$, and the various $\Fpp$, respectively. Moreover, $j_\ell: \SpecOF[1/\ell] \r \SpecOF$ is the open immersion. The remainder of the proof consists in the following steps: we show that the diagram is commutative, that the second row is exact, identify its lower leftmost term and show that the rightmost vertical map is injective. This implies that the left square is cartesian, hence the theorem follows.

We write $\iota$ and $\iota_\ell$ for the open immersions of $U \cap \SpecOF[1/\ell]$ into $\SpecOF[1/\ell]$ and $U$, respectively. 
By \refle{excreal} and the exactness of $j_\ell^*$ we have
$$(j_\ell^* j_{!*} j^* M)_\ell = (\iota_{!*} \iota^* j_\ell^* M)_\ell = \iota_{!*} \iota^* (j_\ell^* M)_\ell.$$
Thus  \refeq{diagramreal} is commutative since every term at the bottom just involves the $\ell$-adic realization of the motive above it, restricted to $\SpecOF[1/\ell]$. 

The exactness of the bottom row is shown separately for each $\ell$, so $\ell$ is fixed for this argument. By the characterization just mentioned, $\iota_{!*} \iota^* (j_\ell^* M)_\ell$ does not change when shrinking $U$, since $j_{!*} j^* M$ is independent of $U$ (as soon as $M$ is smooth over $U$). On the other hand, by the exactness of the $\ell$-adic realization functor (\refax{tstructureRealizations}) $(j_\ell^* M)_\ell$ is a perverse sheaf on $\SpecOF[1/\ell]$, so is a locally constant sheaf (shifted into degree $-1$) on a suitable small open subscheme. Hence we may assume that $\iota^* (j_\ell^* M)_\ell$ is a locally constant sheaf in degree $-1$. By \refsect{perverse}, $ \iota_{!*} \iota^* (j_\ell^* M)_\ell = (\R^0 \iota_* \iota^* (j_\ell^* M)_\ell[-1])[+1]$, so the lower row is the exact cohomology sequence belonging to the distinguished triangle of sheaves on $\SpecOF[1/\ell]$
$$\R^0 \eta_\ell {}_* (M_\eta)_\ell \r \R \eta_\ell {}_* (M_\eta)_\ell \r (\R^1 \eta_\ell {}_* (M_\eta)_\ell)[-1].$$ 
Here $\eta_\ell: \SpecF \r \SpecOF[1/\ell]$ is the generic point. As is well-known, there is an isomorphism
\eqn
\mylabel{eqn_pfiso}
D := \R^1 \eta_\ell {}_* \eta_\ell^* A \stackrel \cong \r \bigoplus_{\pp \nmid \ell} {i_\pp}_* i_\pp^* \R^1 \eta_\ell {}_* \eta_\ell^* A =: \bigoplus B_\pp
\xeqn
for any generically locally constant constructible $\ell$-adic sheaf $A$, such as $M_\ell[-1]$. Indeed, the adjunction map $a: D \r \prod_{\pp \nmid \ell} B_\pp$ factors over the direct sum: note that $(\oplus B_\pp) / \ell^n = \oplus (B_\pp / \ell^n)$ and likewise with the product. Then
$$\Hom(D, \oplus B_\pp) = \varprojlim_n \Hom(D/\ell^n, \oplus (B_\pp / \ell^n)) \subset \varprojlim_n \Hom(D/\ell^n, \prod (B_\pp / \ell^n))$$ 
and to see that $a$ lies in the left hand subgroup, it is enough to consider the $\Z/\ell^n$-sheaves $D / \ell^n$ etc. The corresponding map $\H^1(\Gal(F), A / \ell^n) \r \prod \H^1 (I_\pp, A / \ell^n)$ (Galois cohomology of the inertia groups) factors over the direct sum, since the left hand term agrees with $\H^1(\Gal(F' / F), A)$ for some finite extension $F'/F$. This uses that $A / \ell^n$ is constructible. The extension $F' / F$ is ramified in finitely many places (only), so the claimed factorization follows. This implies \refeq{pfiso} and thus the exactness of the lower row of \refeq{diagramreal}. By \refle{charH1fglobal} and \refle{locconst}, the factors in the lower left-hand term of \refeq{diagramreal} agree with  $\H^1_{f, \backslash \crys}(F, \eta^* M_\ell[-1])$. 

To show that the rightmost vertical map of \refeq{diagramreal} is an injection, let $a = (a_\pp)_{\pp \in \SpecOF}$ be a nonzero element of the rightmost upper term. Only finitely many $a_\pp$ are nonzero. Pick some $\ell$ not lying under any of these prime ideals $\pp$. Then the image of $a$ in $\oplus_{\pp \nmid \ell} \H^0((\pH^0 i_\pp^* \eta_* \eta^* M)_\ell)$ is nonzero by \refle{injectivity}.
\xpf

\lemm \mylabel{lemm_prepinjectivity}
Let $M$ be a mixed motive over $S$ such that $j^* M$ is smooth for some open immersion $j: U \r S$. Then both maps $\H^0 (j_{!*} j^* M) \r \H^0(j_* j^* M) \r \H^0(\eta^*M)$ are injective.
\xlemm
\pf
Indeed the kernels are $\H^{-1}(\pH^0 i^* j_* j^* M) = 0$ and $\oplus_{\pp \in U} \H^0(i_\pp^! M)$, which vanishes since $i_\pp^! M$ sits in cohomological degree $+1$, for $M$ is smooth around $\pp \in U$ (\refle{basic}). 
\xpf

\subsection{$f$-cohomology via $K$-theory of regular models} \mylabel{sect_fcohomologyII}

\defi \mylabel{defi_H1fImage}
Let $X_\eta$ be a smooth and projective variety over $F$. Let $X / \OF$ be any projective model, i.e., $X \x_\OF F = X_\eta$. Then we define
$$\H^i(X_\eta, n)_\OF := \im (\H^i(X, n) \r \H^i(X_\eta, n)).$$ 
\xdefi

Recall that we are assuming the axioms of Sections \ref{sect_axiomaticdescription}, \ref{sect_realizations}, and \ref{sect_conjprop}; the full force of mixed motives will be made use of in the sequel. 

\theo \mylabel{theo_summaryH1f}
The above is well-defined, i.e., independent of the choice of the model $X$. More precisely we have natural isomorphisms:
$$\H^0(\eta_{!*} \h^{i-1}(X_\eta, n)[1]) = \left \{
\begin{array}{cl}
\H^{i}(X_\eta, n)_{\OF} & i<2n \\
\CH^n(X_\eta)_{\Q,\hom} & i=2n 
\end{array} \right.$$
Moreover
$$\H^{-1}(\eta_{!*} \h^{i-1}(X_\eta, n)[1]) = \H^0 (\h^{i-1}(X_\eta, n)).$$
\xtheo

When $X$ is regular, the definition and the statement are due to Beilinson \cite[Lemma 8.3.1]{Beilinson:Notes}. In this case one has $\H^{i}(X_\eta, n)_{\OF} = \im K'_{2j-i}(X)_\Q^{(j)} \r K'_{2j-i}(X_\eta)_\Q^{(j)}$, but that expression does in general depend on the choice of the model \cite{deJeu:Appendix, deJeu:Further}. An extension of Beilinson's definition to all Chow motives over $F$ due to Scholl is discussed in the theorem below. We first provide a preparatory lemma.

\lemm \mylabel{lemm_excimage}
Let $M \in \MM(\SpecOF)$ be a mixed, generically smooth motive with strictly negative weights (\refde{genericallysmooth}). Let $j: U \r \SpecOF$ be an open non-empty immersion such that $M|_U$ is smooth. The natural map $j_{!*} j^* M \r \eta_* \eta^* M$ gives rise to an isomorphism 
$$\H^0(j_{!*} j^* M) = \Im \left ( \H^0(M) \r \H^0(\eta_* \eta^* M) \right).$$ 
\xlemm

\pf
By \refle{prepinjectivity}, $\H^0(j_* j^* M) \r \H^0(\eta_* \eta^* M)$ is injective. Hence it suffices to show $\H^0(j_{!*} j^*M) = \im (\H^0 M \r \H^0 (j_* j^* M))$. Let $i$ be the complement of $j$. From \refeq{local1}, \refeq{local2}, we have a commutative exact diagram 
\scriptsize
$$
\xymatrix{
&
\H^0(j_! j^* M) \ar[r]^\alpha \ar@{->>}[d]
&
\H^0(M) \ar[r] \ar[d] &
\H^0(i_* i^* M) \\
0 = \H^{-1}(i_* \pH^0 i^* j_* j^* M) \ar[r] & 
\H^0(j_{!*}j^* M) \ar@{>->}[r] \ar[d] & 
\H^0(j_* j^* M)
\\
&
\H^{1}(i_* \pH^{-1} i^* j_* j^* M) =0
}
$$
\normalsize
The indicated vanishings are because of $t$-structure reasons and \refax{cohomdim}, respectively. It remains to show that $\alpha$ is surjective. As $i^* M$ is concentrated in cohomological degrees $[-1, 0]$ (\refax{exactness}), there is an exact sequence
$$0 = \H^1(\pH^{-1}i^* M) \r \H^0(i^* M) \r \H^0(\pH^0 i^* M).$$
However $\H^0(\pH^0 i^* M)=0$ as $i^*$ preserves negative weights (\refax{respectweights}) and by strictness of the weight filtration and compatibility with the $t$-structure  (\refax{weight}). 
\xpf

\pf 
Let $j: U \r \SpecOF$ be an open nonempty immersion (which exists by smoothness of $X_\eta$) such that $X_U$ is smooth over $U$. By definition of $\eta_{!*}$ and Lemmas \ref{lemm_eta} and \ref{lemm_smooth}, the left hand term in the theorem agrees with $\H^0(j_{!*} \h^i (X_U, n))$. In the sequel, we write $M := \h^{i}(X, n)$ and $M_\eta := \eta^* [-1] M = \h^{i-1}(X_\eta, n)$.

We first do the case $i < 2n$. The spectral sequences 
$$\H^a(\h^b (X, n)) \Rightarrow \H^{a+b}(X, n), \, \H^a(\h^b (X_\eta, n)) \Rightarrow \H^{a+b}(X_\eta, n)$$
resulting from repeatedly applying truncation functors of the motivic $t$-structure converge since the cohomological dimension is finite (\refax{cohomdim} over $F$, \refle{cohomdimOF} over $\OF$). By \refle{cohomdimOF}, $\H^i(-)$, applied to mixed motives over $\OF$, is non-zero for $i \in \{-1, 0, 1\}$ only. We thus have to consider two exact sequences. The exact functor $\eta^* [-1]$ maps to similar exact sequences for motivic cohomology over $F$ (the indices work out properly, see \refle{eta}):

\eqn \mylabel{eqn_ss1}
\xymatrix{
0 \ar[r] & K \ar[r] \ar[d] & \H^{i}(X,n) \ar[r] \ar[d] & \H^{-1}(\h^{i+1}(X, n)) \ar[d] \ar[r] & 0 \\
0 \ar[r] & K_\eta \ar[r] & \H^{i}(X_\eta, n) \ar[r]  & \H^0(\h^{i}(X_\eta, n)) \stackrel{\ref{axio_weight}, \ref{bsp_weightTate}}= 0 \ar[r]  & 0
}
\xeqn

\eqn \mylabel{eqn_ss2}
\xymatrix{
0 \ar[r] & \H^1(\h^{i-1} (X, n))  \ar[r] \ar[d] & K \ar[r] \ar[d] & \H^0 (M) \ar[d] \ar[r] & 0 \\
0 \ar[r] & \H^2(\h^{i-2}(X_\eta, n)) \stackrel{\ref{axio_cohomdim}}= 0 \ar[r] & K_\eta \ar[r] & \H^1(M_\eta)  \ar[r]  & 0
}
\xeqn
Here, $K$ and $K_\eta$ are certain $E_3$-terms of the spectral sequences above. The rightmost vertical map in \refeq{ss1} is injective as one sees by combining \refeq{localgeneric} with the left-exactness of $i_\pp^!$. Hence 
\eqnarr
\H^{i}(X_\eta, n)_{\OF} & = & \im (\H^{i}(X, n) \r \H^i(X_\eta, n)) = \im (K \r K_\eta) \\
& = & \im (\H^0 (M) \r \H^1(M_\eta))
\xeqnarr
The motive $M = \h^{i}(X, n)$ is a generically smooth (mixed) motive by \refle{smooth}. (Recall that this uses the decomposition axiom \ref{axio_smoothprojective} for smooth projective varieties.) By \refbs{weightTate}, its weights are strictly negative. Thus \refle{excimage} applies and the case $i < 2n$ is shown.

We now do the case $i = 2n$. The motive $j^* M$ is pure of weight zero (\refbs{weightTate}), hence by strictness of the weight filtration for motives over $\OF$ and \refeq{local1}, \refeq{local2} the same is true for $E := j_{!*} j^* M$. (This is an avatar of \cite[Cor. 5.3.2]{BBD}.) Thus $\pH^{1} i^! E$ has strictly positive weights because of \refax{respectweights} and the compatibility of weights and the motivic $t$-structure, i.e., $\wt \pH^{1} (-) \subset \wt (-) + 1$. Therefore $\H^0(\pH^{1} i^! E) = 0$. Here $i$ is any closed immersion. The localization triangle \refeq{localgeneric} yields
$$\H^0 (E) \stackrel{\alpha}\r \H^0(\eta_* \eta^* E) \stackrel{\text{\refeq{morphismsChow}}}= \CH^n(X_\eta)_{\Q, \hom} \r \oplus_\pp \H^1(i_\pp^! E) = \oplus \H^0(\pH^1(i^! E)) = 0.$$ 
Therefore, $\alpha$ is surjective. The injectivity of $\alpha$ is \refle{prepinjectivity}.

To calculate $\H^{-1}(\eta_{!*} M_\eta [1])$, let $j: U \r \SpecOF$ be as above. The natural map $\H^{-1}(\SpecOF, j_{!*} j^* M) \r \H^{-1}(U, j^* M)$ is an isomorphism by the exact cohomology sequence belonging to \refeq{local2}. Thus we have to show 
$$\H^{-1}(\SpecOF, j_* j^* M) = \H^{-1}(\SpecOF, \eta_* \eta^* M).$$
This follows from the localization axiom \ref{axio_localization} and $i_\pp^! M$ being in cohomological degree $+1$ for all points $\pp$ in $U$ (\refle{basic}), so that $\H^0(\Fpp, i_\pp^! M) = \H^{-1}(\Fpp, i_\pp^! M) = 0$.
\xpf

By a theorem of Scholl \cite[Thm. 1.1.6]{Scholl:Integral}, there is a unique functorial and additive (i.e., converting finite disjoint unions into direct sums) way to extend the definition of $\H^i(X_\eta, n)_\OF$ as the image of the $K$-theory of a \emph{regular} proper flat model (\refde{H1fImage}) to all Chow motives over $F$, in particular to ones of smooth projective varieties $X_\eta / F$ that do not possess a regular proper model $X$. The following theorem compares this definition with the one via intermediate extensions. 

\theo \mylabel{theo_Scholl} 
Let $h_\eta$ be a direct summand in the category of Chow motives of $h(X_\eta, n)$ where $X_\eta / F$ is smooth projective.  Let $i \in \Z$ be such that $i -2n <0$. Let $\iota: \PureMot_\rat(F) \r \DMgm(F)$ be the embedding. Then, the group
$$\H^i(h_\eta)_\OF := \H^0 (\eta_{!*} ( \pH^{i-2n-1} (\iota (h_\eta)) [1])).$$
is well-defined and agrees with the aforementioned definition by Scholl.
\xtheo
\pf
Recall $\iota (h(X_\eta, n)) = \M(X_\eta, n)[2n] \in \DMgm(F)$. We first check that the group is well-defined: let $X / \OF$ be a projective model of $X_\eta$. By \refle{spreadout}, there is some model $M \in \MM(\OF)$ of $\pH^{i-2n-1} \iota (h_\eta)[1]$ and an open subscheme $U$ of $\SpecOF$ such that $M$ is a direct summand of $\pH^{i-1} \M(X)(n)$ and such that $X \x U$ is smooth over $U$. Then $\h^{i-1} (X, n)$ is a smooth motive when restricted to $U$ (\refle{smooth}). Hence so is $M$. Thus $\eta_{!*}$ can be applied to $(\pH^{i-2n-1} \iota (h_\eta)) [1]$.

The assignment $h_\eta \mapsto \H^0 (\eta_{!*} (\pH^{i-2n-1} \iota (h_\eta)) [1])$ is functorial and additive and 
$h(X_\eta)(n)$ maps to 
$$\H^0 (\eta_{!*} (\pH^{i-1} \M(X_\eta, n)) [1]) \stackrel{\text{\ref{theo_summaryH1f}}} = \H^i(X_\eta, n)_\OF.$$
Thus the two definitions agree by Scholl's theorem.
\xpf

\bibliography{bib}  

\end{document}

%% file: header_categories.tex
		\newcommand{\category}[1]{\mathbf{#1}}
	
    \newcommand{\D}{\category {D}} 

    \newcommand{\gm}{\mathrm{gm}}

		\newcommand{\Shvv}{\category{Shv}} 
		\newcommand{\Shv}[2]{\Shvv_{#1}\left( #2 \right)} 

    \newcommand{\DM}{\category {DM}}

    \newcommand{\DMgmbeff}{\DM^\eff_\gm} 
    \newcommand{\DMgmeff}{{\DMgmbeff}} 
    \newcommand{\DMeff}{{\DM_\eff}} 
    \newcommand{\DMgm}{\DM_\gm} 
    \newcommand{\DATM}{\category{DATM}} 

    \newcommand{\MM}{\category {MM}} 
    \newcommand{\PureMot}{\category {M}} 


%% file: header_sections.tex
	    \newtheoremstyle{Normal}
  {}
  {}
  {}
  {}
  {\bfseries}
  {.}
  { }
  {}

    \theoremstyle{Normal} 

    \newtheorem{Defi}{Definition}[section]

    \newtheorem{Conj}[Defi]{Conjecture}
    \newtheorem{Bem}[Defi]{Remark}
    \newtheorem{Bsp}[Defi]{Example}
     
    \newtheorem{Ques}[Defi]{Question}

		\theoremstyle{remark}
	
    \theoremstyle{plain} 

    \newtheorem{Satz}[Defi]{Proposition}
    \newtheorem{Theo}[Defi]{Theorem}
    \newtheorem{Folg}[Defi]{Corollary}
    \newtheorem{Lemm}[Defi]{Lemma}
    \newtheorem{Axio}[Defi]{Axiom}

\newcommand{\refsect}[1]{Section \ref{sect_#1}}
\newcommand{\conj}{\begin{Conj}} 			\newcommand{\xconj}{\end{Conj}}											    
\newcommand{\ques}{\begin{Ques}} 			\newcommand{\xques}{\end{Ques}}											    
\newcommand{\axio}{\begin{Axio}} 			\newcommand{\xaxio}{\end{Axio}}											    \newcommand{\refax}[1]{Axiom \ref{axio_#1}}
\newcommand{\bem}{\begin{Bem}} 			\newcommand{\xbem}{\end{Bem}}											    \newcommand{\refbe}[1]{Remark \ref{bem_#1}}
\newcommand{\defi}{\begin{Defi}} 			\newcommand{\xdefi}{\end{Defi}}										\newcommand{\refde}[1]{Definition \ref{defi_#1}}
\newcommand{\lemm}{\begin{Lemm}}			\newcommand{\xlemm}{\end{Lemm}}											\newcommand{\refle}[1]{Lemma \ref{lemm_#1}}
\newcommand{\satz}{\begin{Satz}}			\newcommand{\xsatz}{\end{Satz}}										\newcommand{\refsa}[1]{Proposition \ref{satz_#1}}
\newcommand{\theo}{\begin{Theo}}			\newcommand{\xtheo}{\end{Theo}}											\newcommand{\refth}[1]{Theorem \ref{theo_#1}}
\newcommand{\bsp}{\begin{Bsp}}				\newcommand{\xbsp}{\end{Bsp}}												\newcommand{\refbs}[1]{Example \ref{bsp_#1}}
\newcommand{\folg}{\begin{Folg}}				\newcommand{\xfolg}{\end{Folg}}
				 			
\newcommand{\cor}{\begin{Folg}}				\newcommand{\xcor}{\end{Folg}}
\newcommand{\mycomment}{}

\newcommand{\eqnarra}{\begin{eqnarray}}				\newcommand{\xeqnarra}{\end{eqnarray}}
\newcommand{\eqnarr}{\begin{eqnarray*}}				\newcommand{\xeqnarr}{\end{eqnarray*}}

\newcommand{\eqn}{\begin{equation}} 		\newcommand{\xeqn}{\end{equation}}
\newcommand{\refeq}[1]{(\ref{eqn_#1})}
